\documentclass [12pt]{article}
\usepackage {amsfonts}
\usepackage {mathrsfs}
\usepackage {amsmath}
\usepackage {latexsym}
\usepackage {amssymb}
\usepackage {amsthm}
\usepackage {amscd}
\date{}
\renewcommand{\baselinestretch}{1.2}
\title{\bf Classification of quasi-affine Generalized Dynkin Diagrams with Rank $ 5$ }
\author{ \small Zhengtang Tan $^{a}$,  Shouchuan Zhang $^{b}$   \\
\small $a$.School of Engineering and Design, Hunan Normal University\\ Changsha  410081,   P.R. China \\
\small $b$. Department  of Mathematics,   Hunan University\\
Changsha  410082,  P.R. China\\
\small {\tt Emails:  z9491@sina.cn (SZ);   1843186255@qq.com (ZTT)} }
\date{}

\begin{document}
\newtheorem{Proposition}{Proposition}[section]
\newtheorem{Theorem}[Proposition]{Theorem}
\newtheorem{Definition}[Proposition]{Definition}
\newtheorem{Corollary}[Proposition]{Corollary}
\newtheorem{Lemma}[Proposition]{Lemma}
\newtheorem{Example}[Proposition]{Example}
\newtheorem{Remark}[Proposition]{Remark}

\maketitle


\begin {abstract}  All quasi-affine  connected Generalized Dynkin Diagram with rank $= 5$ are found. All quasi-affine  Nichols (Lie braided) algebras with rank $ 5$  are also found.
\end {abstract}
\section {Introduction and Preliminaries}\label {s0}

Nichols algebras play a fundamental role in the classification of finite-dimensional complex pointed Hopf algebras by means of the lifting method developed by Andruskiewitsch and Schneider \cite {AS02, AS10, AHS08}.
Heckenberger \cite {He06a, He05} classified arithmetic root systems.   Heckenberger \cite {He06b} proved a {\rm  GDD} is  arithmetic  if and only if corresponding matrix is a finite Cartan matrix for {\rm  GDD}s of  Cartan types. W. Wu,   S. Zhang and   Y.-Z. Zhang \cite {WZZ15b} proved a Nichols Lie braided algebra is a finite dimensional  if and only if its {\rm  GDD}, which   fixed
parameter is of finite order, is  arithmetic.

In order to classify finite dimensional Nichols algebras, Heckenberger \cite {He06a, He05} introduced {\rm  GDD}s and  classified arithmetic {\rm  GDD}s.
$\mathfrak B(V)$ is finite dimensional if and only if {\rm  GDD} of $V$ is arithmetic. Zhengtang Tan and   Shouchuan Zhang \cite {TZ24}
introduced quasi-affine {\rm  GDD}s  and  found all quasi-affine  connected Generalized Dynkin Diagrams with rank $> 5$.

We now recall some basic concepts of the graph theory (see \cite {Ha}).
Let $\Gamma _1$ be a non-empty set and $\Gamma _2 \subseteq \{  \{ u, v\} \mid u, v \in \Gamma _1, \hbox { with } u \not= v \} \subseteq 2 ^{\Gamma_1}.$ Then $\Gamma = (\Gamma _1, \Gamma_2)$ is called a graph;  $\Gamma_1$ is called the vertex set of  $\Gamma$; $\Gamma_2$ is called the edge set of  $\Gamma$; Element $\{u, v\} \in \Gamma_2$ is called an edge. Let $F$ be an algebraically closed field of characteristic zero and $F^* := \{x \mid x\in F, x \not= 0 \}$.  If   $\{x_1,   \cdots,   x_n\}$ is  a basis of  vector space  $V$ and
$C(x_i\otimes  x_j) = q_{ij} x_j\otimes x_i$ with $q_{ij} \in F^*$,
then $V$  is called a braided vector space of diagonal type, $\{x_1,   \cdots,   x_n\}$
is called  canonical basis and $(q_{ij})_{n\times n}$  is called braided matrix. Let
$\widetilde{q}_{i j}:= q_{ij}q_{ji}$ for  $i, j \in \{1, 2, \cdots, n \}$ with $i \not=j$. If let $\Gamma _1= \{1, 2, \cdots, n \}$ and $\Gamma _2 = \{  \{ u, v\} \mid \widetilde{q}_{u v} \not= 1,   u, v \in \Gamma _1, \hbox { with } u \not= v \}$, then $\Gamma = (\Gamma _1, \Gamma_2)$ is a graph. Set $q_{ii}$ over vertex $i$ and $\widetilde{q}_{i j}$ over edge $\{i, j\}$ for $i, j \in \Gamma _1$ with $i \not=j$.
Then $\Gamma = (\Gamma _1, \Gamma_2)$ is called a Generalized Dynkin Diagram of braided vector space $V,$ written as {\rm  GDD} in short(see \cite [Def. 1.2.1] {He05}).
If $\Delta (\mathfrak{B}(V))$ is an arithmetic root system.   then we call its {\rm  GDD}  an arithmetic {\rm  GDD}  for convenience.

Let $(q_{ij})_{n \times n }$ be a braided  matrix. If  $q_{ij} q_{ji}  \left \{  \begin{array}{ll}
 \not= 1,  & \mbox {when }  \mid j-i \mid = 1\\
  =1,   & \mbox {when}  \mid j-i \mid \not= 1    \\
\end{array}\right. $ for any $1\le i\not= j \le n,$ then  $(q_{ij})_{n \times n }$ is called a  chain or labelled  chain.
If $(q_{ij})_{n \times n }$ is  a  chain and \begin {eqnarray} \label {ppe1}(q_{11} q_{1, 2} q_{2, 1} -1)(q_{11} +1)=0;
 (q_{n, n } q_{n, n-1} q_{n-1, n} -1)(q_{n,n} +1)=0;\end {eqnarray} i.e.
  \begin {eqnarray} \label {ppe2}
   q_{ii} +1= q_{i, i-1} q_{i-1, i}q_{i, i+1} q_{i+1, i}-1=0
   \end  {eqnarray}
    \begin {eqnarray} \label {ppe3}
   \mbox { or  }
  q_{ii}q_{i, i-1} q_{i-1, i}=q_{ii} q_{i, i+1} q_{i+1, i}=1,\end  {eqnarray}  $1<i < n$,
   then  the braided matrix's {\rm  GDD} is called a simple chain (see \cite [Def.1]{He06a}).  Conditions (\ref {ppe1}), (\ref {ppe2}) and (\ref {ppe3}) are called simple chain conditions.  Let
$q:= q_{n, n} ^2 q_{n -1, n}q_{n, n-1}$.

 For $0\le j \le n$ and $0<i_1 < i_2 < \cdots < i_j \le n,$ let $C_{n,q, i_1, i_2, \cdots, i_j }$ denote  a simple chain which satisfies the condition: $\widetilde{q}_{i, i-1}=q$ if and only if  $i \in \{i_1, i_2, \cdots, i_j\}$,
where
$\widetilde{q} _{1, 0}:= \frac {1}{q_{11}^2\widetilde{q}_{12}}$. For example,
 $n \ge 2$,

$\begin{picture}(100,      15)
\put(27,      1){\makebox(0,     0)[t]{$\bullet$}}
\put(60,      1){\makebox(0,      0)[t]{$\bullet$}}
\put(93,     1){\makebox(0,     0)[t]{$\bullet$}}
\put(159,      1){\makebox(0,      0)[t]{$\bullet$}}
\put(192,     1){\makebox(0,      0)[t]{$\bullet$}}
\put(28,      -1){\line(1,      0){30}}
\put(61,      -1){\line(1,      0){30}}
\put(130,     1){\makebox(0,     0)[t]{$\cdots\cdots\cdots\cdots$}}
\put(160,     -1){\line(1,      0){30}}
\put(22,     -15){1}
\put(58,      -15){2}
\put(91,      -15){3}
\put(157,      -15){$n-1$}
\put(191,      -15){$n$}
\put(22,     10){$q$}
\put(58,      10){$q$}
\put(91,      10){$q$}
\put(157,      10){$q$}
\put(191,      10){$q$}
\put(40,      5){$q^{-1}$}
\put(73,      5){$q^{-1}$}
\put(172,     5){$q^{-1}$}
\put(210,        -1)  {$, q \in F^{*}/\{1, -1\}$. }
\end{picture}$\\

is  $C_{n,q, i_1, i_2, \cdots, i_j }$ with $j=0$.\\ \\

$\begin{picture}(100,      15)
\put(27,      1){\makebox(0,     0)[t]{$\bullet$}}
\put(60,      1){\makebox(0,      0)[t]{$\bullet$}}
\put(93,     1){\makebox(0,     0)[t]{$\bullet$}}
\put(159,      1){\makebox(0,      0)[t]{$\bullet$}}
\put(192,     1){\makebox(0,      0)[t]{$\bullet$}}
\put(28,      -1){\line(1,      0){30}}
\put(61,      -1){\line(1,      0){30}}
\put(130,     1){\makebox(0,     0)[t]{$\cdots\cdots\cdots\cdots$}}
\put(160,     -1){\line(1,      0){30}}
\put(22,     -15){1}
\put(58,      -15){2}
\put(91,      -15){3}
\put(157,      -15){$n-1$}
\put(191,      -15){$n$}
\put(22,     10){$q^{-1}$}
\put(58,      10){$q^{-1}$}
\put(91,      10){$q^{-1}$}
\put(157,      10){$q^{-1}$}
\put(191,      10){$-1$}
\put(40,      5){$q^{}$}
\put(73,      5){$q^{}$}
\put(172,     5){${q}$}
\put(210,        -1)  {$, q \in F^{*}/\{1\}$. }
\end{picture}$\\

is  $C_{n,q, i_1, i_2, \cdots, i_j }$ with $j=n$.

For the convenience, we let
 $C_{1,q, i_1, i_2, \cdots, i_j }$  denote the {\rm  GDD} with length $1$ satisfied the following conditions:   $q= q_{11}$ when $q_{11} \not=-1$; $q$ can be any non-zero number   when $q_{11} =-1$.

Every connected sub{\rm  GDD} of
every   arithmetic {\rm  GDD} in Row 1-10 in Table C is  called a classical {\rm  GDD}.  Every connected     {\rm  GDD} which is not classical is called an exception.

If omitting every vertex in a connected {\rm  GDD} with rank $n>1$  is an arithmetic {\rm  GDD} and this {\rm  GDD} is not an arithmetic  {\rm  GDD}, then this {\rm  GDD} is called  a quasi-affine  {\rm  GDD} over  connected sub{\rm  GDD} of this {\rm  GDD} with rank $n-1$, quasi-affine  {\rm  GDD} in short. In this case, Nichols algebra $\mathfrak B(V)$
and Nichols Lie braded algebra $\mathfrak L(V)$ are said to be quasi-affine.
In other word, if a {\rm  GDD} is quasi-affine  of a braided vector space $V$ which   fixed
parameter is of finite order, then Nichols algebra and Nichols Lie braded algebra
of every proper sub{\rm  GDD} are finite dimensional with $\dim   \mathfrak B { (V)} = \infty$
and $\dim   \mathfrak L { (V)} = \infty$.

In this paper,  using Table $A1$, $A2$, $B$ and $C$ in \cite {He06a, He05}, we find all quasi-affine  connected Generalized Dynkin Diagram with rank $=5$. We also  find all quasi-affine   Nichols  algebras and quasi-affine  Nichols Lie braided algebras with rank $ 5$.

\section {Properties about arithmetic {\rm  GDD} } \label {s2}


\begin {Lemma} \label {mainlemma} A {\rm  GDD} is a classical {\rm  GDD} if and only if it is one of classical Type 1-7, Here
classical types are  listed as follows:\\

 {\ }\ \ \ \ \ \ \ \ \ \ \ \ \ \ \ \ \ \ \ \ \ \ \ \   \ \ \ \ \ \  $\begin{picture}(100,      15)

\put(-125,      -1){ {\rm   Type   1},  $2\le n$.}
\put(80,      1){\makebox(0,      0)[t]{$\bullet$}}

\put(48,      -1){\line(1,      0){33}}
\put(10,      1){\makebox(0,     0)[t]{$C_{n-1, q, i_1, i_2, \cdots, i_j }$}}

\put(-18,     10){$$}
\put(0,      5){$$}
\put(22,     10){$$}
\put(50,      5){$q^{-2}$}

\put(68,      10){$q^2$}

 \ \ \ \ \ \ \ \ \ \ \ \ \ \ \ \ \ \ \
  \ \ \ \ \ \ \ \ \ \ \ \ \ \ \ \ \ \ \ { }$q \in F^{*}\setminus \{1, -1\}$, $0\le j\le n-1.$

 \put(220,      -1) {}
\end{picture}$\\ \\

 {\ }\ \ \ \ \ \ \ \ \ \ \ \ \ \ \ \ \ \ \ \ \ \ \ \   \ \ \ \ \ \  $\begin{picture}(100,      15)

\put(-125,      -1){ {\rm   Type   2}. $2\le n$.}

\put(90,      1){\makebox(0,      0)[t]{$\bullet$}}

\put(58,      -1){\line(1,      0){33}}
\put(27,      1){\makebox(0,     0)[t]{$C_{n-1, q^2, i_1, i_2, \cdots, i_j }$}}

\put(-8,     10){$$}
\put(0,      5){$$}
\put(22,     10){$$}
\put(70,      5){$q^{-2}$}

\put(88,      10){$q$}

  \ \ \ \ \ \ \ \ \ \ \ \ \ \ \ \ \ \ \
  \ \ \ \ \ \ \ \ \ \ \ \ \ \ \ \ \ \ \ $q \in F^{*}\setminus \{1, -1\}$.  $0\le j\le n-1.$
\end{picture}$\\

 {\ }\ \ \ \ \ \ \ \ \ \ \ \ \ \ \ \ \ \ \ \ \ \ \ \   \ \ \ \ \ \  $\begin{picture}(100,      15)

\put(-125,      -1){ {\rm   Type   3}, $2\le n$.}

\put(80,      1){\makebox(0,      0)[t]{$\bullet$}}

\put(48,      -1){\line(1,      0){33}}
\put(27,      1){\makebox(0,     0)[t]{$C_{n-1, q^{-2}, i_1, i_2, \cdots, i_j }$}}

\put(-18,     10){$$}
\put(0,      5){$$}
\put(22,     10){$$}
\put(60,      5){$q^{2}$}

\put(78,      10){$-q^{-1}$}

  \ \ \ \ \ \ \ \ \ \ \ \ \ \ \ \ \ \ \
  \ \ \ \ \ \ \ \ \ \ \ \ \ \ \ \ \ \ \ $q \in F^{*}\setminus \{1, -1\}$, $0\le j\le n-1.$

\end{picture}$\\ \\

 {\ }\ \ \ \ \ \ \ \ \ \ \ \ \ \ \ \ \ \ \ \ \ \ \ \   \ \ \ \ \ \  $\begin{picture}(100,      15)

\put(-125,      -1){{\rm   Type   4}. $2\le n$. }

\put(80,      1){\makebox(0,      0)[t]{$\bullet$}}

\put(48,      -1){\line(1,      0){33}}
\put(27,      1){\makebox(0,     0)[t]{$C_{n-1, -q^{-1}, i_1, i_2, \cdots, i_j }$}}

\put(-18,     10){$$}
\put(0,      5){$$}
\put(22,     10){$$}
\put(60,      5){$-q$}

\put(78,      10){$q$}

  \ \ \ \ \ \ \ \ \ \ \ \ \ \ \ \ \ \ \
  \ \ \ \ \ \ \ \ \ \ \ \ \ \ \ \ \ \ \ $q^3 =1$. $0\le j\le n-1.$

\end{picture}$\\
\\

\ \ \ \     \ \ \ \ \ \  $\begin{picture}(100,      15)

\put(-45,      -1){ {\rm   Type   5}. $3\le n$.}
\put(124,      1){\makebox(0,      0)[t]{$C_{n-2, q^{}, i_1, i_2, \cdots, i_j }$}}
\put(190,     -11){\makebox(0,     0)[t]{$\bullet$}}
\put(190,    15){\makebox(0,     0)[t]{$\bullet$}}
\put(162,    -1){\line(2,     1){27}}
\put(190,      -14){\line(-2,     1){27}}

\put(120,      10){$$}

\put(135,      5){$$}

\put(155,     10){$$}

\put(160,     -20){$q^{-1}$}
\put(165,      15){$q^{-1}$}

\put(193,      -12){$q$}
\put(193,      18){$q$}

 \ \ \ \ \ \ \ \ \ \ \ \ \ \ \ \ \ \ \ \ \ \ \ \ \ \ \ \ \ \ \ \ \ \ \ \ \ \
  \ \ \ \ \ \ \ \ \ \ \ \ \ \ \ \ \ \ \ $q\not= 1,$    $0\le j\le n-2.$

\put(215,        -1)  { }
\end{picture}$\\ \\

\ \ \ \      \ \ \ \ \ \  $\begin{picture}(100,      15)

\put(-45,      -1){{\rm   Type   6}. $3\le n$. }

\put(104,      1){\makebox(0,      0)[t]{$C_{n-2, q^{}, i_1, i_2, \cdots, i_j }$}}
\put(170,     -11){\makebox(0,     0)[t]{$\bullet$}}
\put(170,    15){\makebox(0,     0)[t]{$\bullet$}}

\put(142,    -1){\line(2,     1){27}}
\put(170,      -14){\line(-2,     1){27}}

\put(170,      -14){\line(0,     1){27}}

\put(100,      10){$$}
\put(120,      5){$$}

\put(127,     10){$$}

\put(140,     -20){$q^{-1}$}
\put(145,      15){$q^{-1}$}

\put(178,      -20){$-1$}

\put(178,      0){$q^{2}$}

\put(175,      10){$-1$}

 \ \ \ \ \ \ \ \ \ \ \ \ \ \ \ \ \ \ \ \ \ \ \ \ \ \ \ \ \ \ \ \ \ \ \ \ \ \
  \ \ \ \ \ \ \ \ \ \ \ \ \ \ \ \ \ \ \  $q^2 \not= 1$. $0\le j\le n-2.      $
\put(195,        -1)  {,  }
\end{picture}$\\

{\ }\!\!\!\!\!\!
{\rm   Type   7}. $1\le n$.
$C_{n, q^{-1}, i_1, i_2, \cdots, i_j }.$
$q \not= 1$,  $0\le j\le n$.

\end {Lemma}
{\bf Proof.} We only consider the case $n>4$ since every  {\rm  GDD} classical  type   with $n<5$ is sub{\rm  GDD} of a {\rm  GDD} classical type with  $n> 4.$

Necessity is obviously. Now we show the sufficiency.

 Type   1.   Type   1 is classical  by  {\rm  GDD} $1$ of Row $10$ in Table C when  $2 \le j \le n-1$.   Type   1 is classical  by  {\rm  GDD} $1$ of Row $9$ in Table C when  $j=0$.  Type   1 is classical  by {\rm  GDD} $2$ of Row $9$  in Table C and by {\rm  GDD} $1$ of Row $7$  in Table C
 when  $j=1$.

Type 2.  Type 2 is classical  by  {\rm  GDD} $1$ of Row $4$ in Table C when $q^4\not=1$ and   $1 \le j \le n-1$. Type 2 is classical  by  {\rm  GDD} $1$ of Row $3$ in Table C when $q^2\not=1$ and $j=0$. Type 2 is classical  by  {\rm  GDD} $1$ of Row $3$ in Table C when $1 \le j \le n-1.$

Type 3.  Type 3 is classical  by  {\rm  GDD} $2$ of Row $4$ in Table C when $q^4\not=1$ and  $1 \le j \le n-1$.

Type 4.  Type 4 is classical  by  {\rm  GDD} $1$, $2$ of Row $6$ in Table C when  $1 \le j \le n-1$.  Type 4 is classical  by  {\rm  GDD} $1$ of Row $5$ in Table C when  $j=0$.

Type 5.  Type 5 is classical  by  {\rm  GDD} $3$ of Row $10$ in Table C when  $1 \le j \le n-2$ with $q \not=1$.  Type 5 is classical  by  {\rm  GDD} $1$ of Row $8$  in Table C when  $j=0$ with $q^2\not=1$.

Type 6.  Type 6 is classical  by  {\rm  GDD} $2$ of Row $10$ in Table C when  $1 \le j \le n-2$.  Type 1 is classical  by  {\rm  GDD} $3$ of Row $9$ in Table C when  $j=0$. $\Box$

The sub{\rm  GDD} of the right hand of  Type i in above Lemma is called head type i, written as {\rm h}-type i with $i=1,2,\cdots, 6.$ {\rm h}-type 7 has two forms. Right end is $-1$:    $\begin{picture}(100,      15)  \put(0,      -1){ }

\put(60,      1){\makebox(0,      0)[t]{$\bullet$}}

\put(28,      -1){\line(1,      0){33}}



\put(22,     10){}
\put(38,      5){$q^{-1}$}

\put(60,      10){$-1^{}$}

  \ \ \ \ \ \ \ \ \ \ \ \ \ \ \ \ \ \ \ {  }
\end{picture}$\\
 written $T6$; right end is $-1$:    $\begin{picture}(100,      15)  \put(-20,      -1){ }

\put(60,      1){\makebox(0,      0)[t]{$\bullet$}}

\put(28,      -1){\line(1,      0){33}}



\put(22,     10){}
\put(38,      5){$q^{-1}$}

\put(60,      10){$q^{}$}

  \ \ \ \ \ \ \ \ \ \ \ \ \ \ \ \ \ \ \ { written  T5. }
\end{picture}$\\ The left end of every {\rm  GDD}  in above  Lemma is called the tail of the {\rm  GDD}.

{Remark:}
If $q\in R_3$, then  Type   1 turn into the following:\\

  \ \ \ \ \ \  $\begin{picture}(100,      15)

\put(80,      1){\makebox(0,      0)[t]{$\bullet$}}

\put(48,      -1){\line(1,      0){33}}
\put(27,      1){\makebox(0,     0)[t]{$C_{n-1, q^{}, i_1, i_2, \cdots, i_j }$}}

\put(-18,     10){$$}
\put(0,      5){$$}
\put(22,     10){$$}
\put(60,      5){$q$}

\put(78,      10){$q^{-1}$}

 \ \ \ \ \ \ \ \  \ \ \ \ \ \ \ \ \ \ \ \ \ \ \ \ \ \ \ {$,q^3 =1$.  }

\end{picture}$\\

If $q\in R_3$, then Type   3 turn into the following:\\

  \ \ \ \ \ \  $\begin{picture}(100,      15)

\put(80,      1){\makebox(0,      0)[t]{$\bullet$}}

\put(48,      -1){\line(1,      0){33}}
\put(27,      1){\makebox(0,     0)[t]{$C_{n-1, q^{-1}, i_1, i_2, \cdots, i_j }$}}

\put(-18,     10){$$}
\put(0,      5){$$}
\put(22,     10){$$}
\put(60,      5){$q$}

\put(78,      10){$q$}

 \ \ \ \ \ \ \ \  \ \ \ \ \ \ \ \ \ \ \ \ \ \ \ \ \ \ \ {$,q^3 =1$.  }

\end{picture}$\\

If set   $\xi = -q^{-1}$,
then  Type   3 turn into Type 2: \\

   \ \ \ \ \ \  $\begin{picture}(100,      15)

\put(90,      1){\makebox(0,      0)[t]{$\bullet$}}

\put(58,      -1){\line(1,      0){33}}
\put(27,      1){\makebox(0,     0)[t]{$C_{n-1, \xi ^2, i_1, i_2, \cdots, i_j }$}}

\put(-8,     10){$$}
\put(0,      5){$$}
\put(22,     10){$$}
\put(70,      5){$\xi ^{-2}$}

\put(88,      10){$\xi $}

  \ \ \ \ \ \ \ \ \ \ \ \ \ \ \ \ \ \ \   \ \ \ \ \ \ \ \ \ \ \ \ \ \ \ \ \ \ \ {}$\xi  \in F^{*}\setminus \{1, -1\}$.  $0\le j\le n-1.$
\end{picture}$\\ \\
That is, Type 2 and Type 3 are the same.

If $q\in R_3$, then  Type   3 turn into the following: \\ \\
{\ }  \ \ \ \  \ \ \ \ \ \ \ \ $\begin{picture}(100,      15)

\put(80,      1){\makebox(0,      0)[t]{$\bullet$}}

\put(48,      -1){\line(1,      0){33}}
\put(47,      1){\makebox(0,     0)[t]{$C_{n-1, q^{-1}, i_1, i_2, \cdots, i_j }$}}

\put(-18,     10){$$}
\put(0,      5){$$}
\put(22,     10){$$}
\put(60,      5){$q$}

\put(78,      10){$-q$}

 \ \ \ \ \ \ \ \  \ \ \ \ \ \ \ \ \ \ \ \ \ \ \ \ \ \ \ {$,q^3 =1$.  $\Box$ }

\end{picture}$ \\

A {\rm  GDD} is called a simple cycle  if  it is a cycle and omitting every vertex in the  {\rm  GDD}
is a simple chain. \\

{\ } \ \ \ \ \ \
  \ \ \ \ \ \
  \ \ \ \ \ \
 $\begin{picture}(100,      15)

\put(-98,     10){{\rm  GDD}:}
\put(-28,     10){$\hbox {{\rm h}-Type s} + $}
\put(0,      5){$$}
\put(22,     10){ $C_{u, \xi^{}, i_1, i_2, \cdots, i_j }$}
\put(40,      5){$$}

\put(98,      10){$+ \hbox {{\rm h}-Type t}$}

\put(98,      10){}

\end{picture}$\\
 is called a bi-classical {\rm  GDD}, written as Type s-Type t,  where fixed parameters are  matched as  in
Lemma \ref {mainlemma} with $1\le s, t \le 6.$

A vertex is called an end vertex if its degree  is 1, i.e. there exists only one vertex which connects it.

\begin {Definition} \label {6.1.-2}
{ {\rm (i)}} A connected arithmetic chain   is called a quasi-classical {\rm  GDD}  if omitting an end vertex  is a classical {\rm  GDD} which tail is other end vertex. Further more, the other end vertex is called the tail of the quasi-classical {\rm  GDD}.

{\rm (ii)} A connected  arithmetic  non-chain   is called a quasi-classical {\rm  GDD} if the following conditions hold:  omitting a vertex is a connected classical {\rm  GDD}; omitting an other vertex is also a connected classical {\rm  GDD}; the two tails of the two classical {\rm  GDD}s  are the same. Further more, the two tails of the two classical {\rm  GDD}s is called the tail of the quasi-classical {\rm  GDD}.

{ {\rm (iii)}} A connected arithmetic chain   is called a semi-classical {\rm  GDD}  if omitting an end vertex  is a simple  chain.

{\rm (iv)} A connected  arithmetic  non-chain   is called a semi-classical {\rm  GDD} if the following conditions hold:  omitting a vertex is a connected simple chain; omitting an other vertex is also a connected simple chain.

{\rm (v)} If adding T5 on tail of a  quasi-classical  {\rm  GDD} is an arithmetic {\rm  GDD}   then the  quasi-classical {\rm  GDD} is called continual on the tail via T5.  Otherwise  adding T5 on tail of a  quasi-classical  {\rm  GDD} is called a quasi-affine  continual {\rm  GDD} on the tail via T5.  For T6  we can similarly define these.


\end {Definition}


Remark: Every semi-classical {\rm  GDD} is  a quasi-classical {\rm  GDD}.
 Every  classical {\rm  GDD} is  a semi-classical {\rm  GDD}.

A quasi-classical {\rm  GDD} which is not  semi-classical is called a strict quasi-classical {\rm  GDD}. A semi-classical {\rm  GDD} which is not classical is called a strict semi-classical {\rm  GDD}.

Omitting the longest  simple chain which contains the tail of a quasi-affine {\rm  GDD} is called the head of the {\rm  GDD}.

\begin {Definition} \label {6.1.-3}

{\rm (i)} {\rm h}-Type i  adding on tail of  semi-classical {\rm  GDD} is called classical + semi-classical, $1\le i \le 4.$
{\rm h}-Type 5  adding on tail of  semi-classical {\rm  GDD} is called classical + semi-classical  when  semi-classical {\rm  GDD} is continual via T5 ;  {\rm h}-Type 6  adding on tail of  semi-classical {\rm  GDD} is called classical + semi-classical when  semi-classical {\rm  GDD} is continual via T6.

{\rm (ii)}
The tail of  a semi-classical {\rm  GDD} adding on tail of a semi-classical {\rm  GDD} is called bi-semi-classical, if omitting every vertex of the new  {\rm  GDD} is arithmetic.

 \end {Definition}

 Remark:  Classical + semi-classical can be obtained as follows.
  Adding a vertex on vertex connected tail  such that its tail contained in {\rm h}-Type 5  if the tail of a semi-classical {\rm  GDD} is not $-1$.
   Adding a vertex on vertex connected tail  such that its tail contained in {\rm h}- Type 6  if the tail of a semi-classical {\rm  GDD} is  $-1$.

\begin {Lemma}\label {2.20}
A chain  {\rm  GDD}  with rank $n=4$ is not  arithmetic  if there exist two places  where simple chain conditions do not hold except\\

  {\ }\ \ \ \ \ \ \ \ \ \ \ \  $
$\\ \\
are not arithmetic {\rm  GDD}s containing \ \ \ \ \
  $%
$\\
\\
 is not arithmetic.

\end {Lemma}

\begin {Lemma}\label {2.61}
 If non-classical {\rm  GDD}  with rank $n=4$ \\ \\

 $%
$

\end {Lemma}

\begin {Lemma} \label {2.77}   All strict
semi-classical chain with rank $n=4$ are listed:\\

   {\ }\ \ \ \ \ \ \ \ \ \ \ \  $%
$

 It is the same as {\rm  GDD} $1$ of Row $14$.\\

 {\ }\ \ \ \ \ \ \ \ \ \ \ \  $%
$
\\

{\rm  (r)   }
 classical {\rm  GDD} satisfied conditions.

\end {Lemma}

\begin {Lemma} \label {2.78}
Chain  {\rm  GDD} with rank $n=4$ is not  arithmetic if there exist only one place  where simple chain condition do not hold and $q_{22} \not= -1$, $ \widetilde{q}_{32}q_{22} \not= 1$ or $q_{33} \not= -1$, $ \widetilde{q}_{32}q_{33} \not= 1$,
except\\

  \ \ \ \ \ \  $%
$\\

\end {Lemma}

\begin {Lemma} \label {2.79}
 Non-classical {\rm GDD} with rank $n=4$ is not  arithmetic if there exist two different vertexes  $i$  and $j$ with $q_{ii} \not= -1$
and $\widetilde{q}_{ij} \not= 1 $ such that  $q_{ii} \widetilde{q}_{ij} \not=1$ and $i$ is not a middle vertex (i.e. $i \not= 2,3$) when {\rm GDD} is a chain  except\\

  \ \ \ \ \ \  $%
$\\
\\

\end {Lemma}

\begin {Lemma} \label {2.80}
 A non-classical {\rm GDD} with rank $n=4$ is not  arithmetic  if {\rm GDD} is \\ \\

$%
$\\

\end {Lemma}

\begin {Lemma} \label {2.81}
 A non-classical {\rm GDD} with rank $n=4$ is not  arithmetic if {\rm GDD} is \\ \\

$%
$\\

\end {Lemma}

\begin {Lemma} \label {2.82}
 A non-classical {\rm GDD} with rank $n=4$ is not   arithmetic if {\rm GDD} is \\ \\

$%
$\\
\\

 It is the same as {\rm GDD} $2$ of Row $14$.\\

$%
$\\
\\

\end {Lemma}

\begin {Lemma}\label {2.90}  If {\rm GDD} with rank $n>4$ is  arithmetic non chain,  then omitting some vertex is connected simple chain.

\end {Lemma}

\begin {Lemma}\label {2.63} Assume rank $n >4.$

{\rm (i)}  All  arithmetic {\rm GDD}s are quasi-classical.

{\rm (ii)}  All  classical {\rm GDD}s are semi-classical.

\end {Lemma}

\begin {Lemma} \label {2.91}
If a {\rm GDD} with rank $5$ over a classical {\rm GDD} is quasi-affine, then the  {\rm GDD} is one of following:
{\rm (i)} Bi-classical. {\rm (ii)}  Simple cycle. {\rm (iii)} Quasi-affine  {\rm GDD} over  a non classical {\rm GDD}.
 \\ \\

{\rm (iv)}\ \ \ \ \ \ \ \ $
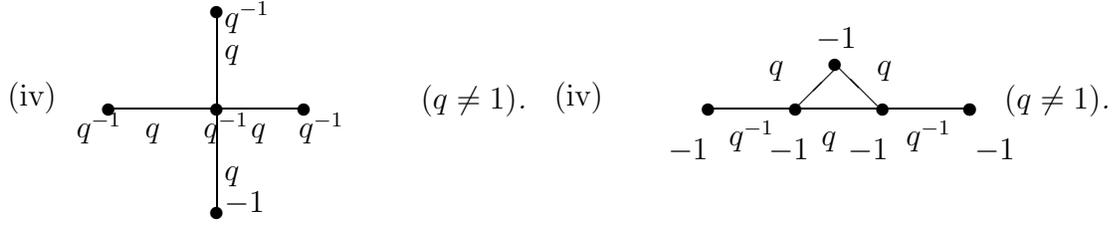
\begin{picture}(100,      15)
\put(-45,      -1){ }
\put(60,      1){\makebox(0,      0)[t]{$\bullet$}}

\put(28,      -1){\line(1,      0){33}}
\put(27,      1){\makebox(0,     0)[t]{$\bullet$}}
\put(-14,      1){\makebox(0,     0)[t]{$\bullet$}}

\put(-14,     -1){\line(1,      0){50}}

\put(-26,     -12){$q^{-1}$}
\put(0,      -12){$q^{}$}
\put(22,     -12){$q^{-1}$}
\put(40,      -12){$q^{}$}
\put(58,      -12){$q^{-1}$}

\put(27,    38){\makebox(0,     0)[t]{$\bullet$}}

\put(27,      0){\line(0,     1){35}}

\put(30,      30){$q^{-1}$}

\put(30,      18){$q^{}$}

\put(27,    -38){\makebox(0,     0)[t]{$\bullet$}}

\put(27,     -38){\line(0,     1){35}}

\put(30,      -40){${-1}$}

\put(30,      -28){$q^{}$}

\put(100,        -1)  { $q \not=1$.   }

\end{picture}$
     \ \ \ \ \ \
  \ \ \ \ \ \
 {\rm (v)}\ \ \
$\begin{picture}(130,      15)

\put(-45,      -1){ }
\put(27,      1){\makebox(0,     0)[t]{$\bullet$}}
\put(60,      1){\makebox(0,      0)[t]{$\bullet$}}
\put(93,     1){\makebox(0,     0)[t]{$\bullet$}}

\put(126,      1){\makebox(0,    0)[t]{$\bullet$}}

\put(126,      1){\makebox(0,    0)[t]{$\bullet$}}
\put(28,      -1){\line(1,      0){33}}
\put(61,      -1){\line(1,      0){30}}
\put(94,     -1){\line(1,      0){30}}

\put(12,    - 20){${-1}$}
\put(35,     - 15){$q^{-1}$}
\put(50,     - 20){${-1}$}
\put(70,     - 15){$q^{}$}

\put(80,      -20){${-1}$}

\put(102,    - 15){$q^{-1}$}

\put(128,     - 20){${-1}$}
\put(135,    - 15){$$}
\put(159,    -  20){$$}

\put(75,    18){\makebox(0,     0)[t]{$\bullet$}}

\put(60,     -1){\line(1,      0){30}}

\put(91,      0){\line(-1,     1){17}}
\put(60,     -1){\line(1,      1){17}}

\put(50,    12){$q^{}$}

\put(68,    22){${-1}$}
\put(91,    12){$q^{}$}

\put(135,        -1) {  $q \in R_3$ . }
\end{picture}$\\ \\

\end {Lemma}

{\bf Proof.}
If it does not hold, then  there exists a quasi-affine    {\rm GDD} over classical Type i which   is  not  bi-classical, simple cycle and over non classical {\rm GDD}s.

{\rm {\rm  (a)   }} If the {\rm GDD}  is a  cycle, then we can obtain a chain such that {\rm h}-Type i is the middle of the chain by omitting a vertex in cycle. This chain is not classical, which is contradiction.

{\rm {\rm  (b)   }} If the {\rm GDD} is a chain,  then the {\rm GDD}  is bi-classical, which is contradiction.

{\rm {\rm  (c)   }} If the {\rm GDD} is not cycle and containing a cycle, then
there is an vertex which does not contain in the cycle. Omitting the vertex in whole is Type 6.
Consequently the {\rm GDD} is a bi-classical or is Case {\rm (v)} since omitting one vertex in head of the Type 6 is classical. It is a contradiction.

{\rm {\rm  (d)   }} If the {\rm GDD} is not chain and does not  contains any cycles, then it is a quasi-affine  over  Type 5.  It has to be  Case {\rm (iv)}.  It is a contradiction. \hfill $\Box$

\begin {Lemma} \label {2.92}

A chain with rank  $n=5$ is not arithmetic if there exist two places  where simple chain conditions do not hold except
 {\rm GDD} $7$ of Row $13$  in Table C. {\rm GDD} $1$ of Row $14$,  {\rm GDD} $1$ of Row $15$.
\end {Lemma}

\begin {Lemma} \label {2.93} Assume rank $n=5.$

{\rm (i)}  Bi-classical {\rm GDD}s are quasi-affine.

{\rm (ii)} Classical  + semi-classical {\rm GDD}s are quasi-affine.

{\rm (iii)} Bi-semi-classical {\rm GDD}s are quasi-affine.

\end {Lemma}
{\bf Proof.} We only need prove that  these {\rm GDD}s are not arithmetic.  If there exists a non chain {\rm GDD} satisfied conditions is  an arithmetic {\rm GDD}, then we obtain a contradiction by Lemma \ref {2.90}. If there exists a chain {\rm GDD} satisfied conditions is  an arithmetic {\rm GDD}, then we obtain a contradiction by Lemma \ref {2.92}.  \hfill $\Box$

\begin {Lemma} \label {2.93'}
 Adding a vertex on tail of a strict quasi-classical {\rm GDD} in Table B is not quasi-affine except\\

$\begin{picture}(130,      15)

\put(-20,      -1){ {\rm  (a)   }}
\put(27,      1){\makebox(0,     0)[t]{$\bullet$}}
\put(60,      1){\makebox(0,      0)[t]{$\bullet$}}
\put(93,     1){\makebox(0,     0)[t]{$\bullet$}}

\put(126,      1){\makebox(0,    0)[t]{$\bullet$}}

\put(159,      1){\makebox(0,    0)[t]{$\bullet$}}
\put(28,      -1){\line(1,      0){33}}
\put(61,      -1){\line(1,      0){30}}
\put(94,     -1){\line(1,      0){30}}
\put(126,      -1){\line(1,      0){33}}

\put(22,     10){$q ^{2}$}
\put(40,      5){$q ^{-2}$}
\put(58,      10){$q ^{2}$}
\put(74,      5){$q ^{-2}$}

\put(91,      10){${q}$}

\put(102,     5){$q ^{-1}$}

\put(120,      10){$-1$}
\put(135,     5){$q ^{3}$}
\put(159,      10){${-3}$}

\put(182,        -1)  { $q^2, q^3 \not=1.$  }
\end{picture}$\\ \\

$\begin{picture}(130,      15)

\put(-20,      -1){ {\rm  (b)   }}
\put(27,      1){\makebox(0,     0)[t]{$\bullet$}}
\put(60,      1){\makebox(0,      0)[t]{$\bullet$}}
\put(93,     1){\makebox(0,     0)[t]{$\bullet$}}

\put(126,      1){\makebox(0,    0)[t]{$\bullet$}}

\put(28,      -1){\line(1,      0){33}}
\put(61,      -1){\line(1,      0){30}}
\put(94,     -1){\line(1,      0){30}}

\put(22,     -20){${-1}$}
\put(40,      -15){${-1}$}
\put(58,      -20){$q ^{}$}
\put(75,      -15){$q ^{-1}$}

\put(87,      10){${-1}$}

\put(102,     5){$q ^{3}$}

\put(128,      10){$q ^{-3}$}

\put(58,    32){\makebox(0,     0)[t]{$\bullet$}}

\put(58,      -1){\line(0,     1){34}}

\put(32,    20){$q^{-1}$}

\put(68,      33){${-1}$}

\put(68,      16){$q ^{-1}$}

\put(28,      -1){\line(1,      1){33}}

\put(166,        -1)  {  $q \in R_4.$  }
\end{picture}$\\  \\ \\

$\begin{picture}(130,      15)

\put(-20,      -1){ {\rm  (c)   }}
\put(27,      1){\makebox(0,     0)[t]{$\bullet$}}
\put(60,      1){\makebox(0,      0)[t]{$\bullet$}}
\put(93,     1){\makebox(0,     0)[t]{$\bullet$}}

\put(126,      1){\makebox(0,    0)[t]{$\bullet$}}

\put(28,      -1){\line(1,      0){33}}
\put(61,      -1){\line(1,      0){30}}
\put(94,     -1){\line(1,      0){30}}

\put(22,     -20){${-1}$}
\put(40,      -15){${-1}$}
\put(58,      -20){${-1}$}
\put(75,      -15){$q ^{-1}$}

\put(91,      10){$q ^{}$}

\put(102,     5){$q ^{-2}$}

\put(128,      10){$q ^{2}$}

\put(58,    32){\makebox(0,     0)[t]{$\bullet$}}

\put(58,      -1){\line(0,     1){34}}

\put(22,    20){$-q ^{-1}$}

\put(68,      33){${-1}$}

\put(68,      16){$q ^{}$}

\put(28,      -1){\line(1,      1){33}}

\put(176,        -1)  { $q \in R_6,$  }
\end{picture}$\\ \\
{\rm  (d)   } Adding {\rm h}-type 7 on tail.

\end {Lemma}
{\bf Proof.}  Let  $\alpha $ be a strict quasi-classical {\rm GDD} in Table B. Let  $\beta$  be a {\rm GDD} by  adding a vertex,   which is not {\rm h}-Type 7, on tail of $\alpha$. Omitting some  vertex  of $\beta$ is Lemma \ref {2.20} {\rm  (a)   }, i.e.  {\rm GDD} $6$ of Row $9$, since $\alpha $ is strict quasi-classical.
Consequently, $\beta$ is {\rm  (a)   } or {\rm  (b)   } or {\rm  (c)   } by  {\rm GDD} $5$ of Row $9$,   {\rm GDD} $4$ of Row $22$ and  {\rm GDD} $2$ of Row $14$.
  \hfill $\Box$

\section {  Main result} \label {s3}

\begin {Lemma} \label {3.1}
All classical + semi-classical  {\rm GDD}s with rank $n=5$ are listed.
\\

$
$

.

\end {Lemma}

\begin {Theorem} \label {Main}
        All connected  quasi-affine {\rm GDD}s with rank $= 5$ are  listed.
{\rm (i)}  Bi-classical {\rm GDD}s.
{\rm (ii)} Simple cycle.
{\rm (iii)}  Quasi-affine  continual {\rm GDD}s.
{\rm (iv)} Classical + semi-classical {\rm GDD}s.
{\rm (v)} Bi-semi-classical {\rm GDD}s.
{\rm (vi)} Other cases
\\

$%
$\\

\end {Theorem}
{\bf Proof.} Considering Lemma \ref {2.91}, we have to find all quasi-affine {\rm GDD}s over
 non classical arithmetic {\rm GDD}s.

 We write the proof according  to the following method.

{\rm  (1)   } If a {\rm GDD} is quasi-affine over two {\rm GDD}s, we write the quasi-affine {\rm GDD} over  behind {\rm GDD}.

 {\rm  (2)   } If a {\rm GDD} is  quasi-affine non chain , we write it in ones over non chain.

{\rm  (3)   } If a non cycle  {\rm GDD} is  quasi-affine with cycle, we write it in ones over {\rm GDD} with cycle.

{\rm  (4)   } If a non cycle  {\rm GDD} is  quasi-affine over classical, we write it in ones over non classical or bi-classical.

{\rm  (5)   }
We do not consider tails of strict quasi-classical {\rm GDD}s by Lemma \ref {2.93'}.

{\rm  (6)   }
We do not consider classical +  strict quasi-classical {\rm GDD}s.
\subsection* {Quasi-affine  {\rm GDD}s  over
{\rm GDD} $1$ of Row $4$} \label {sub3.2} {\ }   {\ }
  {\rm (i)} Omitting  Vertex 1 and adding on  Vertex 4.
 We find all quasi-affine {\rm GDD}s  adding a vertex on Vertex 4 of {\rm GDD} $1$ of Row $4$. We have to consider all  {\rm GDD}s adding a vertex on Vertex 4 of {\rm GDD} $1$ of Row $4$, in which omitting Vertex 1 is an arithmetic {\rm GDD}.  \\

{\rm  (a)   }  $
$\\
It is quasi-affine   since it is not continual.\\

{\rm  (b)   }  $%
$\\  It is {\rm GDD} $21$ of Row $13$ when $q^{}\in R_3.$
It is quasi-affine   since it is not continual when $q^{}\notin R_2 \cup R_3.$
{\rm  (c)   } There are not any other
 cases by Lemma \ref {2.20}.

  (ii) Omitting  Vertex 4 and adding on Vertex 1.\\

$%
$\\
It is quasi-affine   since it is not continual.\\

$%
$\\
It is quasi-affine   since it is not continual.
{\rm  (c)   } There are not any other cases by Lemma \ref {2.79}.

 {\rm (iii)} Cycle. We have to consider the combination between every one in  {\rm (i)} and every one in  {\rm (ii)}.
\\ \\ \\

$\ \ \ \ \  \ \  \ \  \ \ \ \ \ \  %
$
\\ is arithmetic by  {\rm GDD} $5$ of Row $9$.    It is repeated.
\\ \\ \\

$\ \ \ \ \  \ \  \ \  \ \ \ \ \ \  %
$\\ It is the same as case
{\rm  (b)   } to {\rm  (a)   } when $q^{}\in R_4$.\\ \\
Omitting  Vertex 2  {\ }\ \ \ \ \ \ \ \ \ \ \ \  $%
$
\\ is {\rm GDD} $1$ of Row $14$ when $q\in  R_6$;
is not arithmetic  when $q\notin  R_6$ by by Lemma \ref {2.77}.\\ \\
Omitting  Vertex 3  {\ }\ \ \ \ \ \ \ \ \ \ \ \  $%
$
\\ is {\rm GDD} $5$ of Row $14$ when $q\in  R_6$.
  It is repeated. \subsection* {Quasi-affine  {\rm GDD}s  over
{\rm GDD} $1$ of Row $9$} \label {sub3.2}
   {\ }  {\ }  {\rm  (i)   } Omitting Vertex 1 and adding on   Vertex  4.

 {\rm (ii)} Omitting  Vertex  4 and adding on   Vertex  1.

 {\rm (iii)} Cycle.
\\ \\ \\

$\ \ \ \ \  \ \  \ \  \ \ \ \ \   %
$
\\ It is arithmetic by  {\rm GDD} $2$ of Row $9$. It is repeated.
\\ \\ \\

 {\ }

  {\ }

   {\ }

$\ \ \ \ \  \ \  \ \  \ \ \ \ \   %
$\\
It is case  above.
\\ \\ \\

$\ \ \ \ \  \ \  \ \  \ \ \ \ \   %
$
\\ is arithmetic by  Lemma \ref {2.20}{\rm  (b)   } or {\rm GDD} $4$ of Row $14.$
 It is  repeated.
{\rm  (c)   } to {\rm  (b)   }.  It is repeated.
\subsection* {Quasi-affine  {\rm GDD}s  over
{\rm GDD} $2$ of Row $9$} \label {sub3.2} {\ }   {\ }
 {\rm (i)}  Omitting  Vertex  1 and adding on   Vertex  4,

 {\rm (ii)} Omitting  Vertex  4 and adding on   Vertex  1.\\

$%
$\\
 It is repeated. \\

$%
$\\
 It is repeated.
{ {\rm  (j)   }} There are not any other cases by Lemma \ref {2.77}.

 {\rm (iii)} Cycle.\\ \\ \\

$\ \ \ \ \  \ \  \ \  \ \ \ \ \   %
$
\\
is not arithmetic  by Lemma  \ref {2.20}.
\\ \\ \\

$\ \ \ \ \  \ \  \ \  \ \ \ \ \   %
$
\\
is  arithmetic by  {\rm GDD} 2 of  Row $9$. It is quasi-affine.
\\ \\ \\

$\ \ \ \ \  \ \  \ \  \ \ \ \ \   %
$
\\
is not arithmetic  by Lemma  \ref {2.77}  when  $q^4 \not=1$.
{\rm  (a)   } to {\rm  (h)   }. It is repeated. \\ \\   \\

$\ \ \ \ \  \ \  \ \  \ \ \ \ \   %
$
\\
is  arithmetic by  {\rm GDD} 1 of  Row $9$. It is quasi-affine.
{\rm  (b)   } to {\rm  (g)   }.  It is repeated.
\subsection* {Quasi-affine  {\rm GDD}s  over
{\rm GDD} $3$ of Row $9$} \label {sub3.2} {\ }   {\ }
 {\rm (i)} Omitting  Vertex  3 and adding on   Vertex  1\\ \\

$%
$\\ \\  It is repeated. \\ \\

$%
$\\ \\
 It is repeated. \\ \\

$%
$\\ \\
 It is repeated. \\ \\

$%
$\\ \\
 It is repeated.
 { {\rm  (g)   }} There are not any other cases by Lemma \ref {2.20} and  Lemma \ref {2.77}.

 {\rm (ii)}
Omitting  Vertex  1 and adding on   Vertex  3.\\

$%
$\\ \\
Omitting  Vertex  4 is arithmetic. It is quasi-affine  by   Lemma   \ref {2.63}. \\

$%
$\\ \\
 It is repeated.
{\rm  (c)   } There are not any other cases by Lemma \ref {2.60} and Lemma \ref {2.79}.

 {\rm (iii)} Omitting  Vertex  1 and adding on   Vertex  4.\\

$%
$\\ \\
Omitting  Vertex  3 in whole is an arithmetic by Type 1 when $q\in R_4$; by  {\rm GDD} $5$ of Row $9$  when $q\in R_7$;  by  {\rm GDD} $1$ of Row $14$  when $q\in R_6$.
 It is repeated when $q\in  R_7\cup R_6$.\\

$%
$\\ \\
Omitting  Vertex  3 in whole is an arithmetic by  Type 1 when $q\in R_4$;
 is not an arithmetic  by Lemma  \ref {2.77}  when $q\notin R_4$.
{\rm  (c)   } There are not any other cases by Lemma \ref {2.79}.

 {\rm (iv)} Adding on   Vertex  2.\\ \\

$%
$\\
It is empty since Omitting  Vertex  3 is not arithmetic  by Lemma \ref {2.61}.\\ \\

$%
$\\ \\
Omitting  Vertex  3 is not arithmetic Lemma \ref {2.82}.\\ \\

$%
$\\ \\
Omitting  Vertex  1 is  {\rm GDD} $8$ of Row $22$.
 It is repeated.   \\ \\

$%
$\\ \\
Omitting  Vertex  1 is not arithmetic by Lemma \ref {2.59}.\\ \\

$%
$\\ \\
Omitting  Vertex  3 is not arithmetic  Lemma  \ref {2.60}.\\ \\

$%
$\\ \\
Omitting  Vertex  1 is not arithmetic  Lemma  \ref {2.60}.
{\rm  (g)   } There are not any other cases by  Lemma  \ref {2.79}.
\subsection* {Quasi-affine  {\rm GDD}s  over
{\rm GDD} $4$ of Row $9$} \label {sub3.2} {\ }   {\ }
 {\rm (i)} Omitting  Vertex  3 and adding on   Vertex  1. \\ \\

$%
$\\ \\
Omitting  Vertex  4 in whole is an arithmetic. Whole is quasi-affine  by  Lemma   \ref {2.63}.
  \\ \\

$%
$\\ \\
 It is repeated.
{\rm  (c)   } There are not any other cases by Lemma \ref {2.20} and  Lemma \ref {2.77}.

 {\rm (ii)} Adding on   Vertex  2.\\ \\

$%
$\\
Omitting  Vertex  3 is not arithmetic  by  Lemma \ref {2.61} \\ \\

$%
$\\ \\
Omitting  Vertex  3 is not arithmetic by  Lemma \ref {2.60}.
\\ \\

$%
$\\ \\
Omitting  Vertex  3 is not arithmetic by  Lemma \ref {2.82}  when $q\notin R_4$; is not arithmetic by  Lemma \ref {2.82}  when $q\in R_4$ and $q_{55} \not= -1$.
\\ \\

$%
$\\ \\
 Omitting  Vertex  1 in whole is {\rm GDD} 5 of  Row $22$.   It is repeated.
\\ \\

$%
$\\ \\
Omitting  Vertex  1 in whole is not arithmetic  by Lemma \ref {2.60}. \\ \\

$%
$\\ \\
Omitting  Vertex  4 in whole is not arithmetic  by Lemma \ref {2.59}.
{\rm  (g)   } There are not any other cases.

 {\rm (iii)} Omitting  Vertex  3 and adding on   Vertex  4. \\

$%
$\\ \\
Omitting  Vertex  1 in whole is an arithmetic by Type 6 when $q\in R_4$; by  {\rm GDD} $4$ of Row $9$  when $q\in R_7$;  by  {\rm GDD} $2$ of Row $14$  when $q\in R_6$. Omitting  Vertex  1 in whole is not arithmetic when $q\notin R_4\cup R_6 \cup R_7$ by Lemma \ref {2.82}.
 Whole is quasi-affine   when $q\in R_4\cup R_6 \cup R_7$ since it is not continual.  It is repeated in  {\rm GDD} $2$ of Row $14$ when $q\in R_6$.
\\

$%
$\\ \\
Omitting  Vertex  1 is Type 6  when $q\in R_4$.
Omitting  Vertex  1 is not arithmetic  by    Lemma \ref {2.81} when $q\notin R_4$. It is quasi-affine when  $q \in  R_{4}$.
{\rm  (c)   } There are not any other cases by Lemma \ref {2.20}.

 {\rm (iv)} Omitting  Vertex  1 and adding on   Vertex  3.\\

$%
$\\ \\
Omitting  Vertex  4 is simple chain  arithmetic. It is quasi-affine   by Lemma \ref {2.63}. \\

$%
$\\ \\
Omitting  Vertex  4 is Type 1. It is quasi-affine     by Lemma \ref {2.63}.  \\

$%
$\\ \\
  It is repeated.
\subsection* {Quasi-affine  {\rm GDD}s  over
{\rm GDD} $5$ of Row $9$} \label {sub3.2} {\ }   {\ }
 {\rm (i)} Omitting  Vertex  1 and adding on   Vertex  4.\\

$%
$\\
It is quasi-affine  by  Lemma   \ref {2.63}.\\

$%
$\\
 It is repeated. \\

$%
$\\
 It is repeated.
{\rm  (f)   } There are not any other cases by  Lemma \ref {2.20}.

 {\rm (ii)} Omitting  Vertex  4 and adding on   Vertex  1.\\

$%
$\\
 It is repeated.
{\rm  (i)   } There are not any other cases  by Lemma  \ref {2.77}.

 {\rm (iii)}
Cycle.\\ \\ \\

 {\ }

 {\ }

  {\ }

$\ \ \ \ \  \ \  \ \  \ \ \ \ \   %
$
\\
is not arithmetic  by Lemma  \ref {2.77}.
\\ \\ \\

$\ \ \ \ \  \ \  \ \  \ \ \ \ \   %
$
\\
is Type 1. It is quasi-affine.   It is repeated. \\ \\ \\

$\ \ \ \ \  \ \  \ \  \ \ \ \ \   %
$
\\
is not arithmetic by Lemma \ref {2.77}.
{\rm  (b)   } to {\rm  (g)   }.
 It is repeated. \\ \\ \\

$\ \ \ \ \  \ \  \ \  \ \ \ \ \   %
$
\\
is arithmetic by   {\rm GDD} $5$ of Row $9$. It is quasi-affine.
{\rm  (d)   } to {\rm  (e)   }. It is repeated.
{\rm  (e)   } to {\rm  (e)   }. It is repeated.
 {\rm  (e)   } to {\rm  (f)   }. It is repeated.
\subsection* {Quasi-affine  {\rm GDD}s  over
{\rm GDD} $6$ of Row $9$} \label {sub3.2} {\ }   {\ }
 {\rm (i)} Omitting  Vertex  1 and adding on   Vertex  4.\\

$%
$\\
  It is repeated. \\

$%
$\\
It is quasi-affine    since it is not continual.\\

$%
$\\
It is quasi-affine    since it is not continual.\\

$%
$\\
 It is repeated.
  {\rm  (h)   } There are not any other cases by Lemma  \ref {2.20} and Lemma  \ref {2.77}.

 {\rm (ii)} Omitting  Vertex  4 and adding on   Vertex  1.
\\

$%
$\\
It is quasi-affine. \\

$%
$\\  It is repeated.
{\rm  (c)   }   There are not any other cases by Lemma \ref {2.20} and Lemma \ref {2.78}.

 {\rm (iii)} Cycle.\\ \\ \\

$\ \ \ \ \  \ \  \ \  \ \ \ \ \   %
$

{\rm  (d)   } to {\rm  (a)   }. It is repeated. \\ \\ \\

$\ \ \ \ \  \ \  \ \  \ \ \ \ \   %
$
\\ is not arithmetic  by Lemma  \ref {2.20}.
 {\rm  (g)   } to {\rm  (a)   }.  It is repeated.
\subsection* {Quasi-affine  {\rm GDD}s  over
{\rm GDD} $1$ of Row $14$ } \label {sub3.2} {\ }   {\ }
 {\rm (i)} Omitting  Vertex  1 and adding on   Vertex  4.\\

$%
$\\
It is quasi-affine    by  Lemma   \ref {2.63}.\\

$%
$\\
It is quasi-affine    by  Lemma   \ref {2.63}.
{\rm  (c)   } There are not any other cases by Lemma \ref {2.20} and Lemma  \ref {2.77}.

 {\rm (ii)} Omitting  Vertex  4 and adding on   Vertex  1.\\

$%
$\\
 It is repeated.  {\rm  (i)   } There are not any other cases by   Lemma  \ref {2.77}.

 {\rm (iv)} Cycle.\\ \\ \\

$\ \ \ \ \  \ \  \ \  \ \ \ \ \   %
$
\\ is  {\rm GDD} $1$ of Row $14$. It is quasi-affine. \\ \\ \\

$\ \ \ \ \  \ \  \ \  \ \ \ \ \   %
$
\\
 is not arithmetic   by Lemma  \ref {2.78}.
 \\ \\ \\

$\ \ \ \ \  \ \  \ \  \ \ \ \ \   %
$
\\ is not arithmetic    by Lemma  \ref {2.79}.

{\ }\\ \\ \\

$\ \ \ \ \  \ \  \ \  \ \ \ \ \   %
$
\\ is not arithmetic  by  Lemma \ref {2.20}.
 {\rm  (b)   } to {\rm  (f)   }. It is repeated.
\subsection* {Quasi-affine  {\rm GDD}s  over
{\rm GDD} $2$ of Row $14$} \label {sub3.2} {\ }   {\ }
 {\rm (i)} Omitting  Vertex  3 and adding on   Vertex  1.\\ \\

$%
$\\ \\
It is quasi-affine  by  Lemma   \ref {2.63}.

{\rm  (d)   } There are not any other cases by Lemma \ref {2.20} and Lemma  \ref {2.77}.

 {\rm (ii)} Adding on   Vertex  2.\\ \\

$%
$\\ \\
Omitting  Vertex  3 is not arithmetic by Lemma \ref {2.59}. \\ \\

$%
$\\ \\
Omitting  Vertex  4 is not arithmetic by Lemma  \ref {2.59}.
{\rm  (c)   } There are not any other cases by Lemma \ref {2.59}.

 {\rm (iii)} Omitting  Vertex  1 and adding on   Vertex  3.

{\rm  (c)   } There are not any other cases by Lemma \ref {2.59}.

 {\rm (iv)} Omitting  Vertex  1 and adding on   Vertex  4. \\

$%
$\\ \\
Omitting  Vertex  3 in whole is an arithmetic by {\rm GDD} $4$ of Row $14$ Table B. It is repeated. \\

$%
$\\ \\
 It is repeated.\\

$%
$\\ \\
Omitting  Vertex  3 in whole is   an arithmetic by {\rm GDD} $5$ of Row $9$. It is quasi-affine    by Lemma   \ref {2.90}.
{\rm  (d)   } There are not any other cases by  Lemma  \ref {2.59}.
\subsection* {Quasi-affine  {\rm GDD}s  over
{\rm GDD} $3$ of Row $14$} \label {sub3.2} {\ }   {\ }
It is the same as {\rm GDD} $2$ of Row $14$.
\subsection* {Quasi-affine  {\rm GDD}s  over
{\rm GDD} $4$ of Row $14$} \label {sub3.2} {\ }   {\ }
  {\rm (i)} Omitting  Vertex  1 and adding on   Vertex  4.\\

$%
$\\
It is quasi-affine    since it is not continual.
{\rm  (d)   } There are not any other cases by Lemma \ref {2.20} and Lemma  \ref {2.77}.

 {\rm (ii)} Omitting  Vertex  4 and adding on   Vertex  1.\\

$%
$\\
It is quasi-affine     since it is not continual.
{\rm  (d)   } There are not any other cases by Lemma \ref {2.20} and Lemma  \ref {2.77}.

 {\rm (iii)} Cycle.\\ \\ \\

$\ \ \ \ \  \ \  \ \  \ \ \ \ \   %
$
\\
is a simple chain. It is quasi-affine.
\\ \\ \\

$\ \ \ \ \  \ \  \ \  \ \ \ \ \   %
$
\\
is arithmetic by  {\rm GDD} $3$ of Row $22.$
  It is repeated. \\ \\ \\

$\ \ \ \ \  \ \  \ \  \ \ \ \ \   %
$
\\
is not  arithmetic by Lemma \ref {2.20}.\\ \\ \\

$\ \ \ \ \  \ \  \ \  \ \ \ \ \   %
$
\\
is not  arithmetic by Lemma \ref {2.20}.
\subsection* {Quasi-affine  {\rm GDD}s  over
{\rm GDD} $5$ of Row $14$} \label {sub3.2} {\ }   {\ }
 It is the same as {\rm GDD} $1$ of Row $14$.
\subsection* {Quasi-affine  {\rm GDD}s  over
{\rm GDD} $1$ of Row $17$} \label {sub3.2} {\ }   {\ }
 {\rm (i)} Omitting  Vertex  1 and adding on   Vertex  4.
{\rm  (a)   } There are not any  cases by Lemma \ref {2.20} and  Lemma \ref {2.78}.

 {\rm (ii)} Omitting  Vertex  4 and adding on   Vertex  1.\\

$%
$\\
It is quasi-affine   by  Lemma   \ref {2.63}.
{\rm  (i)   } There are not any other cases by  Lemma  \ref {2.77}.

 {\rm (iii)} Cycle. It is  empty.
\subsection* {Quasi-affine  {\rm GDD}s  over
{\rm GDD} $2$ of Row $17$} \label {sub3.2} {\ }   {\ }
 {\rm (i)} Omitting  Vertex  1 and adding on   Vertex  4.
{\rm  (a)   } There are not any other cases by Lemma \ref {2.20}.

 {\rm (ii)} Omitting  Vertex  4 and adding on   Vertex  1.
{\rm  (h)   } There are not any other cases by Lemma \ref {2.77}.

 {\rm (iii)} Cycle. It is  empty.
\subsection* {Quasi-affine  {\rm GDD}s  over
{\rm GDD} $3$ of Row $17$} \label {sub3.2} {\ }   {\ }
 {\rm (i)} Omitting  Vertex  1 and adding on   Vertex  3.
{\rm  (a)   } There are not any  cases by Lemma \ref {2.79}.

 {\rm (ii)} Omitting  Vertex  1 and adding on   Vertex  4.
{\rm  (a)   } There are not any  cases by  Lemma \ref {2.79}.

 {\rm (iii)} Omitting  Vertex  3 and adding on   Vertex  1.

 {\rm (iv)} Adding on   Vertex  2.\\ \\

$
$\\ \\
Omitting  Vertex  3 is not arithmetic  Lemma \ref {2.61}.\\ \\

$%
$\\  \\
Omitting  Vertex  3 is not arithmetic by Lemma \ref {2.82}.
{\rm  (c)   } There are not any other cases.
\subsection* {Quasi-affine  {\rm GDD}s  over
{\rm GDD} $4$ of Row $17$} \label {sub3.2} {\ }   {\ }
 {\rm (i)} Omitting  Vertex  1 and adding on   Vertex  3.

 {\rm (ii)} Omitting  Vertex  1 and adding on   Vertex  4.\\

$%
$\\
 is not arithmetic  by Lemma  \ref {2.20}.
{\rm  (b)   } There are not any other cases by  Lemma \ref {2.82}.

 {\rm (iii)} Omitting  Vertex  3 and adding on   Vertex  1.
{\rm  (a)   } There are not any  cases by   Lemma  \ref {2.77} and  Lemma  \ref {2.20}.

 {\rm (iv)} Adding on   Vertex  2.\\ \\

$%
$\\ \\
Omitting  Vertex  3 is not arithmetic by  Lemma \ref {2.61}.\\ \\

$%
$\\ \\
Omitting  Vertex  3 is not arithmetic by Lemma \ref {2.60}.\\ \\

$%
$\\ \\
Omitting  Vertex  1 is not arithmetic by Lemma \ref {2.82}.
{\rm  (d)   } There are not any other cases.
\subsection* {Quasi-affine  {\rm GDD}s  over
{\rm GDD} $5$ of Row $17$} \label {sub3.2} {\ }   {\ }
 {\rm (i)} Omitting  Vertex  1 and adding on   Vertex  4.

 {\rm (ii)} Omitting  Vertex  4 and adding on   Vertex  1.
{\rm  (a)   } There are not any other cases by Lemma \ref {2.20} and    Lemma  \ref {2.78}.

 {\rm (iii)} Cycle.
\subsection* {Quasi-affine  {\rm GDD}s  over
{\rm GDD} $6$ of Row $17$
} \label {sub3.2} {\ }   {\ }
 {\rm (i)} Omitting  Vertex  1 and adding on   Vertex  4.\\

$%
$\\
It is quasi-affine    by  Lemma   \ref {2.63}.
{\rm  (b)   } There are not any other cases   by Lemma  \ref {2.77} and Lemma  \ref {2.20}.

 {\rm (ii)} Omitting  Vertex  4 and adding on   Vertex  1.
{\rm  (h)   } There are not any other cases by Lemma \ref {2.77}.
 {\rm (iii)} Cycle.
 \\ \\ \\

$\ \ \ \ \  \ \  \ \  \ \ \ \ \   %
$
\\
It is quasi-affine.
\subsection* {Quasi-affine  {\rm GDD}s  over
{\rm GDD} $1$ of Row $18$} \label {sub3.2} {\ }   {\ }

 {\rm (i)} Omitting  Vertex  1 and adding on   Vertex  4.

 {\rm (ii)} Omitting  Vertex  4 and adding on   Vertex  1.

 {\rm (iii)} Cycle.\\ \\ \\

$\ \ \ \ \  \ \  \ \  \ \ \ \ \   %
$
\\ is arithmetic by Type 7.
It is quasi-affine.
\\ \\ \\

$\ \ \ \ \  \ \  \ \  \ \ \ \ \   %
$
\\
is arithmetic by {\rm GDD} $2$ of Row $18$.
It is quasi-affine.
\subsection* {Quasi-affine  {\rm GDD}s  over
{\rm GDD} $2$ of Row $18$
 } \label {sub3.2} {\ }   {\ }
 {\rm (i)} Omitting  Vertex  1 and adding on   Vertex  4.\\
 {\rm (ii)} Omitting  Vertex  4 and adding on   Vertex  1.\\

 $%
$\\
 It is repeated.
{\rm  (h)   } There are not any other cases by   Lemma  \ref {2.77}.

 {\rm (iii)} Cycle.

 {\ } \\ \\ \\

$\ \ \ \ \  \ \  \ \  \ \ \ \ \   %
$
\\
is  arithmetic by {\rm GDD} $1$ of Row $18$. It is quasi-affine.
 \\ \\ \\

$\ \ \ \ \  \ \  \ \  \ \ \ \ \   %
$

{\rm  (b)   } to {\rm  (f)   }.  It is repeated.
\subsection* {Quasi-affine  {\rm GDD}s  over
{\rm GDD} $3$ of Row $18$
} \label {sub3.2} {\ }   {\ }
 {\rm (i)} Omitting  Vertex  4 and adding on   Vertex  1. \\  \\

$%
$\\ \\
 It is repeated.
 {\rm  (g)   } There are not any other cases by  Lemma  \ref {2.77}.

 {\rm (ii)} Omitting  Vertex  1 and adding on   Vertex  4.

 {\rm (iii)} Omitting  Vertex  1 and adding on   Vertex  3.\\

$%
$\\ \\
Omitting  Vertex  4 is arithmetic.
It is quasi-affine    by  Lemma   \ref {2.63}.\\

$%
$\\ \\
  It is repeated.

{\rm  (c)   } There are not any other cases by  Lemma \ref {2.79} and  Lemma  \ref {2.60}.

 {\rm (iv)} Adding on   Vertex  2.\\ \\

$%
$\\ \\
Omitting  Vertex  1 is arithmetic by {\rm GDD} $8$ of Row $20$.
 It is repeated. \\ \\

$%
$\\ \\
 Omitting  Vertex  1 is not arithmetic by Lemma  \ref {2.60}.
\\ \\

$%
$\\ \\ Omitting  Vertex  1 is not arithmetic by Lemma  \ref {2.60}.\\
\\

$%
$\\ \\
 Omitting  Vertex  1 is not arithmetic by Lemma  \ref {2.79}.
\\ \\

$%
$\\ \\
 Omitting  Vertex  4 is
  {\rm GDD} $4$ of Row $18$.   It is repeated.
 \\ \\

$%
$\\ \\
Omitting  Vertex  4 is {\rm GDD} $4$ of Row $18$. It is repeated.
{\rm  (i)   } There are not any other cases  by Lemma \ref {2.79}.
\subsection* {Quasi-affine  {\rm GDD}s  over
{\rm GDD} $4$ of Row $18$} \label {sub3.2} {\ }   {\ }
 {\rm (i)} Omitting  Vertex  1 and adding on   Vertex  3.\\ \\

$%
$\\ \\
 It is repeated.

{\rm  (g)   } There are not any other cases by   Lemma  \ref {2.77}.

 {\rm (ii)} Omitting  Vertex  3 and adding on   Vertex  4.

 {\rm (iii)} Omitting  Vertex  3 and adding on   Vertex  1.\\ \\

$%
$\\ \\
It is quasi-affine   by  Lemma   \ref {2.63}.\\ \\

$%
$\\ \\
  It is repeated.

{\rm  (c)   } there are not any other cases by Lemma \ref {2.20}.

 {\rm (iv)} Omitting  Vertex  1 and adding on   Vertex  2\\ \\

$%
$\\ \\
Omitting  Vertex  3 is  {\rm GDD} $9$ of Row $20$.  It is repeated.
\\ \\

$%
$\\ \\
It is  quasi-affine   by Lemma   \ref {2.90}.
{\rm  (e)   } There are not any other cases by Lemma \ref {2.79}.
\subsection* {Quasi-affine  {\rm GDD}s  over
{\rm GDD} $5$ of Row $18$}\label {sub3.2} {\ }   {\ }
 {\rm (i)} Omitting  Vertex  3 and adding on   Vertex  4.

 {\rm (ii)} Omitting  Vertex  4 and adding on   Vertex  3.

 {\rm (iii)} Omitting  Vertex  3 and adding on   Vertex  1.\\ \\

$%
$\\ \\
It is quasi-affine    by Lemma   \ref {2.90}. \\ \\

$%
$\\

 It is repeated.

{\rm  (c)   } There are not any other cases by Lemma \ref {2.20}.

 {\rm (iv)} Adding on   Vertex  2.\\ \\

$%
$\\ \\
Omitting  Vertex  1 is {\rm GDD} $6$ of Row $18$. It is repeated.
{\rm  (c)   } There are not any other cases  by Lemma  \ref {2.79}.
\subsection* {Quasi-affine  {\rm GDD}s  over
{\rm GDD} $6$ of Row $18$}\label {sub3.2} {\ }   {\ }
 {\rm (i)} Omitting  Vertex  1 and adding on   Vertex  3.\\ \\

$%
$\\ \\
  It is repeated.

{\rm  (g)   }
 There are not any other cases    by Lemma  \ref {2.77}.

 {\rm (ii)} Omitting  Vertex  3 and adding on   Vertex  1.
\\ \\

$%
$\\ \\
 It is repeated.

{\rm  (g)   } There are not any other cases  by Lemma  \ref {2.77}.

 {\rm (iii)} Adding on   Vertex  2\\ \\

$%
$\\ \\
Omitting  Vertex  3 is  {\rm GDD} $10$ of Row $20$.    It is repeated. \\ \\

 {\ }

$%
$\\ \\
{\rm  (e)   } There are not any other cases  by  Lemma \ref {2.79}.

 {\rm (iv)} Omitting  Vertex  1 and adding on   Vertex  4 are the same as  Omitting  Vertex  1 and adding on   Vertex  3.
\subsection* {Quasi-affine  {\rm GDD}s  over
{\rm GDD} $1$ of Row $20$} \label {sub3.2} {\ }   {\ }
 {\rm (i)}  Omitting  Vertex  1 and  adding on   Vertex  4.\\

$%
$\\
 It is repeated.
\\

$%
$\\
 It is repeated.
{\rm  (c)   } There are not any other cases by Lemma  \ref {2.20}.

 {\rm (ii)} Omitting  Vertex  4 and adding on   Vertex  1\\

$%
$\\
 It is repeated.
{\rm  (h)   } There are not any other cases  by Lemma  \ref {2.77}.

 {\rm (iii)} Cycle.
{\rm  (a)   } to {\rm  (d)   }. It is repeated. {\rm  (a)   } to {\rm  (f)   }.  It is repeated.
{\rm  (b)   } to {\rm  (a)   }.  It is repeated.
{\rm  (b)   } to {\rm  (b)   }.  It is repeated.
\subsection* {Quasi-affine  {\rm GDD}s  over
{\rm GDD} $2$ of Row $20$} \label {sub3.2} {\ }   {\ }
 {\rm (i)} Omitting  Vertex  1 and adding on   Vertex  4.\\

$%
$\\
It is quasi-affine  by Lemma   \ref {2.90}.\\

$%
$\\
 It is repeated.

{\rm  (c)   } There are not any other cases by  Lemma \ref {2.20}.

 {\rm (ii)} Omitting  Vertex  4 and adding on   Vertex  1.

 {\rm (iii)} Cycle.
\\ \\

$\ \ \ \ \  \ \  \ \  \ \ \ \ \   %
$

{\rm  (b)   } to {\rm  (b)   }.  It is repeated.
\subsection* {Quasi-affine  {\rm GDD}s  over
{\rm GDD} $3$ of Row $20$} \label {sub3.2} {\ }   {\ }
 {\rm (i)} Omitting  Vertex  1 and adding on   Vertex  4.\\

$%
$\\
It is quasi-affine   by Lemma \ref {2.63}. \\

$%
$\\
 It is repeated.

{\rm  (c)   } There are not any other cases by  Lemma  \ref {2.20}.

 {\rm (ii)}  Omitting  Vertex  4 and adding on   Vertex  1.\\

$%
$\\
It is quasi-affine    by Lemma   \ref {2.90}.

{\rm  (h)   } There are not any other cases by Lemma \ref {2.77}.

 {\rm (iii)} Cycle.
\\ \\ \\

$\ \ \ \ \  \ \  \ \  \ \ \ \ \   %
$
\\ It
is not arithmetic  by Lemma  \ref {2.20}.
{\rm  (b)   } to {\rm  (a)   }. It is repeated.
 {\rm  (b)   } to {\rm  (b)   }. It is repeated.
\subsection* {Quasi-affine  {\rm GDD}s  over
{\rm GDD} $4$ of Row $20$} \label {sub3.2} {\ }   {\ }
 {\rm (i)} Omitting  Vertex  1 and adding on   Vertex  4.

 {\rm (ii)}  Omitting  Vertex  4 and adding on   Vertex  1.\\

$%
$\\
It is quasi-affine   by  Lemma   \ref {2.63}.

{\rm  (h)   } There are not any other cases by   Lemma  \ref {2.77}.

 {\rm (iii)} Cycle.
\\ \\ \\

$\ \ \ \ \  \ \  \ \  \ \ \ \ \   %
$
\\
Omitting  Vertex  2
is simple chain.\\
Omitting  Vertex  3
 {\ }\ \ \ \ \ \ \ \ \ \ \ \  $%
$
\\
is arithmetic by Type 1.
It is quasi-affine. \\ \\ \\

$\ \ \ \ \  \ \  \ \  \ \ \ \ \   %
$
\\
is  {\rm GDD} $6$ of Row $21$ or  Lemma  \ref {2.77}{\rm  (n)   }.
 It is repeated.
\subsection* {Quasi-affine  {\rm GDD}s  over
{\rm GDD} $5$ of Row $20$} \label {sub3.2} {\ }   {\ }
 {\rm (i)} Omitting  Vertex  1 and adding on   Vertex  4.\\

$%
$\\
 It is repeated.
\\

$%
$\\
It is quasi-affine    by  Lemma   \ref {2.63}.

{\rm  (c)   } There are not any other cases  by Lemma  \ref {2.20}.

 {\rm (ii)} Omidding 4 and adding on   Vertex  1.

 {\rm (iii)}  Cycle.

{\rm  (a)   } to {\rm  (a)   }.
It is repeated.\\ \\ \\

$\ \ \ \ \  \ \  \ \  \ \ \ \ \   %
$
\subsection* {Quasi-affine  {\rm GDD}s  over
{\rm GDD} $6$ of Row $20$} \label {sub3.2} {\ }   {\ }
 {\rm (i)} Omitting  Vertex  1 and adding on   Vertex  4.

 {\rm (ii)} Omitting  Vertex  4 and adding on   Vertex  1.

 {\rm (iii)} Cycle.
\\ \\

$\ \ \ \ \  \ \  \ \  \ \ \ \ \   %
$
\subsection* {Quasi-affine  {\rm GDD}s  over
{\rm GDD} $7$ of Row $20$} \label {sub3.2} {\ }   {\ }
 {\rm (i)} Omitting  Vertex  1 and adding on   Vertex  3.
{\rm  (a)   } There are not any  cases by Lemma \ref {2.80}.

 {\rm (ii)} Omitting  Vertex  1 and adding on   Vertex  4.

 {\rm (iii)} Omitting  Vertex  3 and adding on   Vertex  1.\\ \\

$%
$\\ \\
Omitting  Vertex  4 is arithmetic by {\rm GDD} $4$ of Row $20.$
It is quasi-affine  by Lemma   \ref {2.90}.\\ \\

$%
$\\ \\
Omitting  Vertex  4 is not arithmetic  by Lemma \ref {2.77}.
{\rm  (j)   } There are not any other  cases   by Lemma  \ref {2.77}.

 {\rm (iv)} Adding on   Vertex  2.\\ \\

$%
$\\ \\
Omitting  Vertex  1 is  {\rm GDD} $7$ of Row $21$.
 It is repeated.\\ \\

$%
$\\ \\
 Omitting  Vertex  1 is {\rm GDD} $7$ of Row $21$  in Table B.
 It is repeated.\\ \\

$%
$\\ \\
If $\widetilde{q} _{15}\not= q^2$, then
Omitting  Vertex  3 is not arithmetic by non classical.
{\rm  (f)   } There are not any other cases by Lemma \ref {2.79}.
\subsection* {Quasi-affine  {\rm GDD}s  over
{\rm GDD} $8$ of Row $20$} \label {sub3.2} {\ }   {\ }
 {\rm (i)} Omitting  Vertex  1 and adding on   Vertex  3.\\

$%
$\\ \\
Omitting  Vertex  4 in whole is simple chain. It is quasi-affine    by  Lemma   \ref {2.63}.
\\

$%
$\\ \\
Omitting  Vertex  4 in whole is simple chain. It is quasi-affine.
{\rm  (c)   } There are not any other cases by Lemma \ref {2.79} and  Lemma \ref {2.60}.

 {\rm (ii)}  Omitting  Vertex  1 and adding on   Vertex  4.

 {\rm (iii)} Omitting  Vertex  3 and adding on   Vertex  1.\\ \\

$%
$\\ \\
 Omitting  Vertex  4 is  {\rm GDD} $3$ of Row $20.$
 It is repeated.

{\rm  (g)   } There are not any other cases by Lemma  \ref {2.77}.

 {\rm (iv)} Adding on   Vertex  2.\\ \\

$%
$\\ \\
Omitting  Vertex  1 is not arithmetic  Lemma \ref {2.60} when $ \widetilde{q}_{52} \not= q$.\\ \\

$%
$\\ \\
Omitting  Vertex  1 is not arithmetic   when $ \widetilde{q}_{52} q_{55} \not= 1$ and  $q_{55} \not=-1$ Lemma \ref {2.79}.\\ \\

$%
$\\ \\
 Omitting  Vertex  3 is {\rm GDD} $9$ of Row $20.$
 It is repeated. \\ \\

$%
$\\ \\
Omitting  Vertex  1 is {\rm GDD} $3$ of Row $18.$ Omitting  Vertex  4 is {\rm GDD} $8$ of Row $20.$
Omitting  Vertex  3 is {\rm GDD} $3$ of Row $18.$ It is quasi-affine. \\ \\

$%
$\\ \\
Omitting  Vertex   4 is not not arithmetic by  Lemma \ref {2.81}.
{\rm  (h)   } There are not any other cases.
\subsection* {Quasi-affine  {\rm GDD}s  over
{\rm GDD} $9$ of Row $20$} \label {sub3.2} {\ }   {\ }

 {\rm (i)} Omitting  Vertex  1 and adding on   Vertex  4.
 {\rm (ii)} Omitting  Vertex  1 and adding on   Vertex  3. \\ \\

$%
$\\ \\
It is quasi-affine  by  Lemma   \ref {2.63}. \\ \\

$%
$\\ \\
It is quasi-affine   by  Lemma   \ref {2.63}.
{\rm  (c)   } There are not any other cases by   Lemma  \ref {2.20}.

 {\rm (iii)} Adding on   Vertex  2. \\ \\

$%
$\\ \\
  It is quasi-affine  since it is not  continual.\\ \\

$%
$\\ \\
Omitting  Vertex  1 is not arithmetic by Lemma \ref {2.82}. \\ \\

$%
$\\ \\
Omitting  Vertex  1 is arithmetic  by {\rm GDD} $4$ of Row $18$.
Omitting  Vertex  3 is arithmetic by  {\rm GDD} $4$ of Row $18$.
 Omitting  Vertex  4 is arithmetic by  {\rm GDD} $9$ of Row $20$. It is quasi-affine. \\ \\

$%
$\\ \\
Omitting  Vertex  1 is not arithmetic by  Lemma \ref {2.82}.
{\rm  (e)   } There are not any other cases by Lemma \ref {2.79}.

 {\rm (iv)}  Omitting  Vertex  3 and adding on   Vertex  1.\\ \\

$%
$\\ \\
 It is repeated.

{\rm  (g)   } There are not any other cases by Lemma \ref {2.77}.
\subsection* {   Quasi-affine  {\rm GDD}s  over       {\rm GDD} $10$ of Row $20$}\label {sub3.2}
 {\ }   {\ }
 {\rm (i)} Omitting  Vertex  4 and adding on   Vertex  3.\\ \\

$%
$\\ \\
Omitting  Vertex  1 is  {\rm GDD} $1$ of Row $18$. It is quasi-affine.
{\rm  (b)   }  There are not any  cases by Lemma \ref {2.20} and   Lemma \ref {2.78}.

 {\rm (ii)}  Omitting  Vertex  4 and adding on   Vertex  1.

 {\rm (iii)} Adding on   Vertex  2\\ \\

 {\ }

$%
$\\ \\
Omitting  Vertex  1 is arithmetic  by {\rm GDD} $5$ of Row $18$.
Omitting  Vertex  3 is {\rm GDD} $6$ of Row $18$. Omitting  Vertex  4 is {\rm GDD} $10$ of Row $20$.  It is quasi-affine. \\ \\

$%
$\\ \\
Omitting  Vertex  1 is not arithmetic by Lemma  \ref {2.79}.
{\rm  (d)   } There are not any other cases   by Lemma  \ref {2.79}.

 {\rm (iv)} Omitting  Vertex  1 and  adding on   Vertex  4.
\subsection* {Quasi-affine  {\rm GDD}s  over
{\rm GDD} $1$ of Row $21$} \label {sub3.2} {\ }   {\ }
 {\rm (i)} Omitting  Vertex  1 and adding on   Vertex  4.
{\rm  (a)   } There are not any  cases by  Lemma \ref {2.78} and  Lemma \ref {2.20}.
 {\rm (ii)} Omitting  Vertex  4 and adding on   Vertex  1.\\

$%
$\\
 It is quasi-affine   by  Lemma   \ref {2.63}.\\

$%
$\\
It is quasi-affine   by  Lemma   \ref {2.63}.\\

$%
$\\ It is
quasi-affine   by  Lemma   \ref {2.63}.

{\rm  (k)   } There are not any other cases by Lemma  \ref {2.77}.
\subsection* {Quasi-affine  {\rm GDD}s  over
{\rm GDD} $2$ of Row $21$} \label {sub3.2} {\ }   {\ }
 {\rm (i)} Omitting  Vertex  1 and adding on   Vertex  4.
{\rm  (a)   } There are not any  cases by Lemma \ref {2.20}.

 {\rm (ii)} Omitting  Vertex  4 and adding on   Vertex  1.
\subsection* {
{\rm GDD} $3$ of Row $21$} \label {sub3.2}
 {\ }   {\ }   {\rm (i)} Omitting  Vertex  1 and adding on   Vertex  4.\\

$%
$\\
It is quasi-affine  by  Lemma   \ref {2.63}.
{\rm  (c)   } There are not any other cases by Lemma \ref {2.20}.

 {\rm (ii)} Omitting  Vertex  4 and adding on   Vertex  1.\\

$%
$\\
 It is repeated.
{\rm  (h)   } There are not any other cases   by Lemma  \ref {2.77}.

 {\rm (iii)} Cycle.
 \\ \\ \\

$\ \ \ \ \  \ \  \ \  \ \ \ \ \   %
$\\ \\
Omitting  Vertex  2 is not arithmetic    by Lemma  \ref {2.79}.
{\rm  (a)   } to {\rm  (f)   }. It is repeated.\\ \\ \\

$\ \ \ \ \  \ \  \ \  \ \ \ \ \   %
$
 \subsection* {    Quasi-affine  {\rm GDD}s  over
{\rm GDD} $4$ of Row $21$} \label {sub3.2} {\ }   {\ }
 {\rm (i)} Omitting  Vertex  1 and adding on   Vertex  4.\\

$%
$\\
It is quasi-affine    by  Lemma   \ref {2.63}.
{\rm  (c)   } There are not any other cases by Lemma \ref {2.20}.

  {\rm (ii)}   Omitting  Vertex  4 and adding on   Vertex  1.

  {\rm (iii)} Cycle.
\\ \\ \\

$\ \ \ \ \  \ \  \ \  \ \ \ \ \   %
$
\\
Omitting  Vertex  2, 3 are arithmetic. It is quasi-affine.
 \subsection* {    Quasi-affine  {\rm GDD}s  over
{\rm GDD} $5$ of Row $21$} \label {sub3.2} {\ }   {\ }
 {\rm (i)} Omitting  Vertex  1 and adding on   Vertex  4.\\

$%
$\\
It is quasi-affine    by  Lemma   \ref {2.63}.
{\rm  (c)   } There are not any other cases  by Lemma  \ref {2.20}.
 {\rm (ii)} Omitting  Vertex  4 and adding on   Vertex  1.

 {\rm (iii)} Cycle.
\\ \\ \\

$\ \ \ \ \  \ \  \ \  \ \ \ \ \   %
$
\\
Omitting  Vertex  2, 3 are arithmetic. It is quasi-affine.
 \subsection* {    Quasi-affine  {\rm GDD}s  over
{\rm GDD} $6$ of Row $21$} \label {sub3.2} {\ }   {\ }
 {\rm (i)} Omitting  Vertex  1 and adding on   Vertex  4.\\

$%
$\\
It is quasi-affine    by  Lemma   \ref {2.63}.
{\rm  (c)   } There are not any other cases by Lemma \ref {2.20}.

 {\rm (ii)} Omitting  Vertex  4 and adding on   Vertex  1.\\

$%
$\\
It is quasi-affine      by  Lemma   \ref {2.63}.
{\rm  (h)   } There are not any other cases by Lemma \ref {2.77}.

 {\rm (iii)} Cycle.
 \\ \\ \\

$\ \ \ \ \  \ \  \ \  \ \ \ \ \   %
$
\\ It is quasi-affine.
\\ \\ \\

$\ \ \ \ \  \ \  \ \  \ \ \ \ \   %
$
\\
is arithmetic by {\rm GDD} $5$ of Row $20$.
It is quasi-affine.
   \\ \\ \\

$\ \ \ \ \  \ \  \ \  \ \ \ \ \   %
$
\\
is arithmetic by Type 7.
It is quasi-affine.
 \subsection* {    Quasi-affine  {\rm GDD}s  over
{\rm GDD} $7$ of Row $21$} \label {sub3.2} {\ }   {\ }
 {\rm (i)} Omitting  Vertex  1 and adding on   Vertex  3.
{\rm  (a)   } There are not any  cases by
 Lemma \ref {2.80}.
 {\rm (ii)} Omitting  Vertex  1 and adding on   Vertex  4.\\

$%
$\\ \\
 Omitting  Vertex  3 is a simple chain. It is quasi-affine     by  Lemma   \ref {2.63}.\\

$%
$\\ \\
 Omitting  Vertex  3 is a simple chain. It is quasi-affine  by  Lemma   \ref {2.63}.
{\rm  (c)   } There are not any other cases Lemma \ref {2.79}.

 {\rm (iii)} Omitting  Vertex  3 and adding on   Vertex  1.
{\rm  (f)   } There are not any other cases by Lemma  \ref {2.77}.

 {\rm (iv)} Adding on   Vertex  2.\\ \\

$%
$\\ \\
Omitting  Vertex  1 is {\rm GDD} $7$ of Row $20$.
Omitting  Vertex  3 is  Type 5.   Omitting  Vertex  4 is arithmetic by   {\rm GDD} $6$ of Row $18$.  It is quasi-affine. \\ \\

$%
$\\ \\
Omitting  Vertex  1 is {\rm GDD} $7$ of Row $20$.
Omitting  Vertex  3 is   {\rm GDD} $7$ of Row $21$.   Omitting  Vertex  4 is arithmetic by   {\rm GDD} $7$ of Row $20$.  It is quasi-affine. \\ \\

$%
$\\ \\
Omitting  Vertex  3 is not   arithmetic by   Lemma \ref {2.80}.
{\rm  (f)   } There are not any other cases by
 Lemma \ref {2.79}.
 \subsection* {    Quasi-affine  {\rm GDD}s  over
{\rm GDD} $1$ of Row $22$} \label {sub3.2} {\ }   {\ }
 {\rm (i)} Omitting  Vertex  1 and adding on   Vertex  4.
{\rm  (a)   } There are not any  cases by Lemma \ref {2.20} and  Lemma \ref {2.78}.

 {\rm (ii)} Omitting  Vertex  4 and adding on   Vertex  1.\\

$%
$\\
It is quasi-affine    since it is not continuale.

{\rm  (h)   } There are not any other cases by  Lemma  \ref {2.77}.
 \subsection* {    Quasi-affine  {\rm GDD}s  over
{\rm GDD} $2$ of Row $22$} \label {sub3.2} {\ }   {\ }
 {\rm (i)} Omitting  Vertex  1 and adding on   Vertex  4.
{\rm  (c)   } There are not any other cases by  Lemma  \ref {2.20} and   Lemma  \ref {2.77}.

 {\rm (ii)} Omitting  Vertex  4 and adding on   Vertex  1.\\

$%
$\\
It is quasi-affine.
{\rm  (h)   } There are not any other cases  by Lemma  \ref {2.77}.
 \subsection* {    Quasi-affine  {\rm GDD}s  over
{\rm GDD} $3$ of Row $22$} \label {sub3.2} {\ }   {\ }
 {\rm (i)} Omitting  Vertex  1 and adding on   Vertex  4.\\

$%
$\\
It is quasi-affine     by  Lemma   \ref {2.63}.
\\

$%
$\\
It is quasi-affine     by  Lemma   \ref {2.63}.
{\rm  (d)   } There are not any other cases by Lemma  \ref {2.20}.

 {\rm (ii)} Omitting  Vertex  4 and adding on   Vertex  1.\\

$%
$
\\ is  arithmetic   {\rm GDD} $1$ of Row $14$. It is quasi-affine.

 {\ }\\ \\ \\

 {\ }

$\ \ \ \ \  \ \  \ \  \ \ \ \ \   %
$
\\ It is  arithmetic by {\rm GDD} $4$ of Row $14$. It is quasi-affine.
 \\ \\ \\

$\ \ \ \ \  \ \  \ \  \ \ \ \ \   %
$
\\ is  simple chain. It is quasi-affine.
 \subsection* {    Quasi-affine  {\rm GDD}s  over
{\rm GDD} $4$ of Row $22$} \label {sub3.2} {\ }   {\ }
 {\rm (i)} Omitting  Vertex  1 and adding on   Vertex  3.
{\rm  (a)   } There are not any  cases   by Lemma  \ref {2.79}.
 {\rm (ii)} Omitting  Vertex  1 and adding on   Vertex  4.
{\rm  (a)   } There are not any  cases   by Lemma  \ref {2.79}.

 {\rm (iii)} Omitting  Vertex  3 and adding on   Vertex  1.\\ \\

$%
$\\ \\
It is quasi-affine    by Lemma   \ref {2.90}.
{\rm  (d)   } There are not any other cases by Lemma \ref {2.20}.

 {\rm (iv)} Omitting  Vertex  1 and adding on   Vertex  2.\\ \\

$%
$\\ \\
Omitting  Vertex  3 is not arithmetic by   Lemma  \ref {2.79}. {\rm  (b)   } There are not any  cases  by Lemma  \ref {2.79}.
 \subsection* {    Quasi-affine  {\rm GDD}s  over
{\rm GDD} $5$ of Row $22$} \label {sub3.2} {\ }
  {\ } {\rm (i)} Omitting  Vertex  1 and adding on   Vertex  3.

 {\rm (ii)} Omitting  Vertex  1 and adding on   Vertex  4.\\

$%
$\\ \\
 Omitting  Vertex  3  is not arithmetic by Lemma  \ref {2.20}\\

$%
$\\ \\
 Omitting  Vertex  3 is  {\rm GDD} $2$ of Row $9$ or  Lemma  \ref {2.77}{\rm  (a)   }.
It is quasi-affine     by Lemma   \ref {2.90}. \\

$%
$\\ \\
 Omitting  Vertex  3 is {\rm GDD} $7$ of Row $22.$ It is repeated.
{\rm  (d)   } There are not any other cases by Lemma \ref {2.59}.

 {\rm (iii)} Omitting  Vertex  3 and adding on   Vertex  1.

 {\rm (iv)} Omitting  Vertex  2.\\ \\

$%
$\\ \\
Omitting  Vertex  1 is arithmetic Lemma \ref {2.59}{\rm  (b)   } or  {\rm GDD} $4$ of Row $9$.
Omitting  Vertex  3 is the same as Omitting  Vertex  1.
Omitting  Vertex  4  is arithmetic by  {\rm GDD} $5$ of Row $22$. It is quasi-affine.
{\rm  (c)   } There are not any other cases by Lemma \ref {2.59}.
 \subsection* {    Quasi-affine  {\rm GDD}s  over
{\rm GDD} $6$ of Row $22$} \label {sub3.2} {\ }   {\ }
 {\rm (i)} Omitting  Vertex  1 and adding on   Vertex  4.
{\rm  (a)   } There are not any cases   by Lemma  \ref {2.20}.
 {\rm (ii)} Omitting  Vertex  4 and adding on   Vertex  1.\\

$%
$\\
 It is quasi-affine   since it is not  continual.
{\rm  (d)   } There are not any other cases by Lemma \ref {2.20}.
 \subsection* {    Quasi-affine  {\rm GDD}s  over
{\rm GDD} $7$ of Row $22$} \label {sub3.2} {\ }   {\ }
 {\rm (i)} Omitting  Vertex  1 and adding on   Vertex  4.

  {\rm (ii)} Omitting  Vertex  4 and adding on   Vertex  1.

 {\rm (iii)} Cycle.
\\ \\ \\

$\ \ \ \ \  \ \  \ \  \ \ \ \ \   %
$
 \subsection* {    Quasi-affine  {\rm GDD}s  over
{\rm GDD} $8$ of Row $22$} \label {sub3.2}  {\ }  {\ }
 {\rm (i)} Omitting  Vertex  1 and adding on   Vertex  3.

 {\rm (ii)} Omitting  Vertex  1 and adding on   Vertex  4.\\

$%
$\\ \\
Omitting  Vertex  3 is arithmetic by Lemma \ref {2.20}{\rm  (b)   } or  {\rm GDD} $4$ of Row $14.$
It is quasi-affine     by Lemma   \ref {2.90}.\\

$%
$\\ \\
Omitting  Vertex  3 is   {\rm GDD} $5$ of Row $9.$
It is quasi-affine     by Lemma   \ref {2.90}.
{\rm  (c)   } There are not any other cases by Lemma \ref {2.59}.

 {\rm (iii)} Omitting  Vertex  3 and adding on   Vertex  1.

 {\rm (iv)} Adding on   Vertex  2. \\ \\

$%
$\\ \\
Omitting  Vertex  3 is arithmetic by {\rm GDD} $3$ of Row $9$  in Table B or   Lemma \ref {2.59}{\rm  (a)   }.
Omitting  Vertex  4 is arithmetic by {\rm GDD} $8$ of Row $22$  in Table B. Omitting  Vertex  1 is arithmetic by {\rm GDD} $3$ of Row $9$  in Table B. It is quasi-affine.
{\rm  (c)   } There are not any other cases by Lemma \ref {2.59}.

\section {Appendix}

Given a generalized Cartan matrix $A$ we can obtain a Kac-Moody Lie algebra $ \mathfrak g {  (A)   }$. If $A$ is a Cartan matrix, then  $ \mathfrak g {  (A)   }$ is a semi-simple Lie algebra.   If $A$ is a quasi-affine  generalized Cartan matrix, then  $ \mathfrak g {  (A)   }$ is a quasi-affine  Lie algebra, which is infinite dimensional. Given a braided matrix $Q$,  its braided vector $V$,  we can obtain a Nichols  algebra $\mathfrak B(V)$ and  Nichols Lie braided algebra $\mathfrak L (V)$. If the {\rm GDD} is an arithmetic {\rm GDD}  which   fixed
parameter is of finite order, then $\mathfrak B(V)$ and $\mathfrak L (V)$  are  finite dimensional.   If the {\rm GDD} is  quasi-affine   then $\mathfrak B(V)$ and  $\mathfrak L (V)$  are infinite dimensional (see \cite {WZZ15a, WZZ15b, WWZZ}).

Nichols Lie super algebra $\mathfrak L (V)$ and Nichols Lie colour
  algebra $\mathfrak L (V)$ are trivial when $V$  is of  diagonal type since the braiding of $V$ is symmetric.

If a {\rm GDD} is obtained  by relabelling   other {\rm GDD}, then we view the two {\rm GDD} are the same since their YD modules are isomorphic.
For example, \\

 $\begin{picture}(100,      15)

\put(60,      1){\makebox(0,      0)[t]{$\bullet$}}

\put(28,      -1){\line(1,      0){33}}
\put(27,      1){\makebox(0,     0)[t]{$\bullet$}}

\put(-14,      1){\makebox(0,     0)[t]{$\bullet$}}

\put(-14,     -1){\line(1,      0){50}}

\put(-18,     10){$q$}
\put(0,      5){$q ^{-1}$}
\put(22,     10){$q$}
\put(40,      5){$q ^{-1}$}

\put(58,      10){$q ^{-3}$}

\put(122,      -1){and }

  \ \ \ \ \ \ \ \ \ \ \ \ \ \ \ \ \ \ \ {$ q \in R_9$ }
\end{picture}$\ \ \ \ \ \ \ \ \ \ \ \
  \ \ \ \ \ \
 $\begin{picture}(100,      15)

\put(60,      1){\makebox(0,      0)[t]{$\bullet$}}

\put(28,      -1){\line(1,      0){33}}
\put(27,      1){\makebox(0,     0)[t]{$\bullet$}}

\put(-14,      1){\makebox(0,     0)[t]{$\bullet$}}

\put(-14,     -1){\line(1,      0){50}}

\put(-18,     10){ $q ^{-3}$}
\put(0,      5){$q ^{-1}$}
\put(22,     10){$q$}
\put(40,      5){$q ^{-1}$}

\put(58,      10){$q$}

  \ \ \ \ \ \ \ \ \ \ \ \ \ \ \ \ \ \ \ {$ q \in R_9$\ \ \  are the same. }
\end{picture}$
 \subsection {   About semi-classical, continual  and tail  } {\ }   {\ }
Vertex 1 is tail of all {\rm GDD}s with rank 4 except  {\rm GDD} $4$ of Row $9$,  {\rm GDD} $6$ of Row $9$, {\rm GDD} $4$ of Row $17$, {\rm GDD} $5$ of Row $17$.   {\rm GDD} $4$ of Row $18$,    {\rm GDD} $5$ of Row $18$.

{\rm GDD} $1$ of Row $4$ is continual on Vertex  4: {\rm GDD} $21$ of Row $13$ in  Table C via T6. Vertex 4 is a tail .

 {\rm GDD} $1$ of Row $9$ is continual on   Vertex  1: {\rm GDD} $7$ of Row $15$  in  Table C via T5 , $q \in R_5$.   Vertex   4 is a tail .

 {\rm GDD} $2$ of Row $9$ is continual on   Vertex  1: {\rm GDD} $6$ of Row $15$  in  Table C via T6 , $q \in R_5$.   Vertex  4 is a tail .  It is semi-classical.

 {\rm GDD} $3$ of Row $9$ is continual on  Vertex   1: {\rm GDD} $3$ of Row $15$ in  Table C  via T5 , $q \in R_5$.   Vertex  4 is a tail when  $q \in R_4$.

 {\rm GDD} $4$ of Row $9$
  is continual on   Vertex  1: {\rm GDD} $2$ of Row $15$ in  Table C  on   Vertex  3 via T5, $q \in R_5$.   Vertex  3 is a tail  when $q \in  R_{5}$.
  Vertex  4 is a tail  when $q \in  R_{4}$.   Vertex  1 is not a tail .

 {\rm GDD} $5$ of Row $9$
 is semi-classical.  Vertex   4 is a tail  when $q\in R_4$.

  {\rm GDD} $6$ of Row $9$ is continual: {\rm GDD} $1$ of Row $15$  in  Table C via T5  on  Vertex   4, $q\in R_5$.   Vertex  1 is not a tail .

 {\rm GDD} $1$ of Row $14$ is continual: {\rm GDD} $4$ of Row $14$ in  Table C  via T5 , $q\in R_4$.  It is semi- classical.

 {\rm GDD} $2$ of Row $14$ is continual: {\rm GDD} $2$ of Row $14$  in  Table C via T5  on   Vertex  3, $q\in R_4$;  {\rm GDD} $3$ of Row $14$  in  Table C  via  T5  on   Vertex  1,  $q\in R_4$.     Vertex   1 is a tail  when $q\in R_4$.  Vertex   3 is a tail  when $q\in R_4$.

  {\rm GDD} $3$ of Row $14$ is the same as {\rm GDD} $2$ of Row $14$.  Vertex   1 is a tail  when $q\in R_4$.

  {\rm GDD} $4$ of Row $14$ is continual: {\rm GDD} $1$ of Row $14$  in  Table C  via  T5  on   Vertex  1, $q\in R_4$.
  Vertex  4 is a tail  when $q\in R_4$.

  {\rm GDD} $5$ of Row $14$ is the same as {\rm GDD} $1$ of Row $14$.

   {\rm GDD} $1$ of Row $17$ is semi-classical on tail    Vertex  1.

   {\rm GDD} $2$ of Row $17$ is semi-classical on tail   Vertex   1.

    {\rm GDD} $4$ of Row $17$.
  Vertex  1 is not a tail .
  Vertex  3 is a tail .

    {\rm GDD} $5$ of Row $17$.   Vertex  1 is not a tail .
  Vertex  4 is a tail .

    {\rm GDD} $6$ of Row $17$ is semi-classical on tail    Vertex  1.

      {\rm GDD} $1$ of Row $18$
 is continual:
{\rm GDD} $1$ of Row $13$  in  Table C  via  T6  on   Vertex  1; {\rm GDD} $1$ of Row $12$  in  Table C  via  T5  on   Vertex  1; {\rm GDD} $15$ of Row $12$  in  Table C  via  T5  on   Vertex  4; {\rm GDD} $20$ of Row $13$  in  Table C  via  T6  on   Vertex  4.

       {\rm GDD} $2$ of Row $18$ is continual: {\rm GDD} $2$ of Row $12$  in  Table C  via  T5  on   Vertex  1;  {\rm GDD} $3$ of Row $13$  in  Table C  via  T6  on   Vertex  1;  {\rm GDD} $12$ of Row $12$  in  Table C  via  T6  on   Vertex  4. It is semi-classical.

       {\rm GDD} $3$ of Row $18$
 is continual: {\rm GDD} $3$ of Row $12$  in  Table C  via  T5  on   Vertex  1;  {\rm GDD} $6$ of Row $13$  in  Table C  via  T6  on  Vertex   1;  {\rm GDD} $9$ of Row $12$  in  Table C  via  T6  on   Vertex  4; {\rm GDD} $18$ of Row $13$  in  Table C  via  T5  on   Vertex  4. It is semi-classical.

       {\rm GDD} $4$ of Row $18$
   is continual: {\rm GDD} $4$ of Row $12$  in  Table C  via  T5  on   Vertex  3;  {\rm GDD} $10$ of Row $13$  in  Table C  via  T6  on   Vertex  3;  {\rm GDD} $6$ of Row $12$  in  Table C  via  T6  on   Vertex  4; {\rm GDD} $17$ of Row $13$  in  Table C  via  T5  on   Vertex  4. It is semi-classical on   Vertex  3.   Vertex  1 is not a tail .

     {\rm GDD} $5$ of Row $18$
It is continual: {\rm GDD} $15$ of Row $13$  in  Table C  via  T6  on   Vertex  4;  {\rm GDD} $5$ of Row $12$  in  Table C  via  T5  on   Vertex  4.  Vertex   1 is not a tail .

      {\rm GDD} $6$ of Row $18$ is continual: {\rm GDD} $14$ of Row $13$  in  Table C  via  T5  on  Vertex   4;  {\rm GDD} $8$ of Row $12$  in  Table C  via  T6  on  Vertex   4. It is semi-classical on  Vertex   1,  Vertex   3,  Vertex   4, respectively.

            {\rm GDD} $1$ of Row $20$ is continual: {\rm GDD} $11$ of Row $13$  in  Table C  via  T5  on   Vertex  1;  It is semi-classical on tail    Vertex  1.

      {\rm GDD} $2$ of Row $20$ is continual: {\rm GDD} $7$ of Row $13$  in  Table C  via  T5  on  Vertex   1.

       {\rm GDD} $3$ of Row $20$ is continual: {\rm GDD} $8$ of Row $13$  in  Table C  via  T5  on  Vertex   1. It is semi-classical.

       {\rm GDD} $4$ of Row $20$ is continual: {\rm GDD} $5$ of Row $13$  in  Table C  via  T6  on   Vertex  1. 4 is a tail . It is semi-classical on tail  1.

       {\rm GDD} $6$ of Row $20$ is continual: {\rm GDD} $2$ of Row $13$  in  Table C  via  T6  on   Vertex  1.  4 is a tail .

        {\rm GDD} $7$ of Row $20$ is continual: {\rm GDD} $12$ of Row $13$  in  Table C  via  T5  on   Vertex  1;  It is semi-classical.

      {\rm GDD} $8$ of Row $20$ is continual: {\rm GDD} $9$ of Row $13$  in  Table C  via  T6  on   Vertex  1.   Vertex  4 is a tail . It is semi-classical on tail  1.

       {\rm GDD} $9$ of Row $20$ is continual:   {\rm GDD} $13$ of Row $13$  in  Table C  via  T6  on   Vertex  1;  Vertex   4 is a tail . It is semi-classical on tail    Vertex  1.

        {\rm GDD} $10$ of Row $20$ is continual: {\rm GDD} $16$ of Row $13$  in  Table C  via  T6  on   Vertex  1.   Vertex   4 is a tail.

       {\rm GDD} $1$ of Row $21$ is semi-classical on tail   Vertex   1.
       {\rm GDD} $3$ of Row $21$  is semi-classical on tail    Vertex  1.
 {\rm GDD} $6$ of Row $21$  is semi-classical on tail    Vertex  1.
 {\rm GDD} $7$ of Row $21$ is semi-classical on tail   Vertex   1.
   {\rm GDD} $1$ of Row $22$ is semi-classical on tail   Vertex   1.
      {\rm GDD} $2$ of Row $22$ is semi-classical on tail   Vertex   1.
       {\rm GDD} $3$ of Row $22$ is semi-classical on tail   Vertex   1.
             {\rm GDD} $5$ of Row $22$.   Vertex   3 is a tail .
               {\rm GDD} $7$ of Row $22$.   Vertex  4 is a tail .      {\rm GDD} $8$ of Row $22$.  Vertex   3 is a tail .
 \subsection {   About proof of main result} {\ }   {\ }
In order to
  find all  quasi-affine  {\rm GDD}s over an  arithmetic {\rm GDD}  $\alpha $ in Table B, we have to  find all  quasi-affine  {\rm GDD}s which are {\rm GDD}s adding   a vertex on Vertex $1, 2, \cdots, n$, respectively, of  $\alpha $. For example,    we find all  quasi-affine  {\rm GDD}s which are {\rm GDD}s adding   a vertex on Vertex $1$. Choice $i$, $2\le i \le n.$
Let  $\beta $  the {\rm GDD} omitting    Vertex i  of  $\alpha $. We find all  {\rm GDD}s  by adding  a vertex on Vertex 1 of  $\beta $  such that they are  arithmetic. Next we put Vertex 1 into these {\rm GDD}s. Finally we determine if they are quasi-affine.


\end {document}